\documentclass[11pt,english]{article}

\usepackage{anysize}
\marginsize{2cm}{2cm}{2.5cm}{2.4cm}

\usepackage[english]{babel}
\usepackage{amsfonts}
\usepackage{latexsym}
\usepackage{amsmath}
\usepackage{mathrsfs}
\usepackage{amssymb}
 \usepackage{graphpap}
 \usepackage{graphicx}
 \usepackage{epsfig}

\let\phi=\varphi
\newcommand{\Z}{{\mathbb Z}}
\newcommand{\R}{{\mathbb R}}
\newcommand{\N}{{\mathbb N}}
\newcommand{\eps}{\varepsilon}
\newcommand{\IP}{{\mathbb P}}

\newcommand{\IB}{{\mathbb B}}

\newcommand{\G}{{\mathcal G}}
\newcommand{\Zc}{{\mathcal Z}}
\newcommand{\Xc}{{\mathcal X}}
\newcommand{\Dc}{{\mathcal D}}

\newcommand{\Rc}{{\mathcal R}}

\newcommand{\Po}{{\mathtt P}_{\omega}}

\newcommand{\Eo}{{\mathtt E}_{\omega}}

\newcommand{\sig}{\sigma}
\newcommand{\alf}{\alpha}
\newcommand{\gam}{\gamma}

\newcommand{\sgn}{\operatorname{sgn}}

\newcommand{\qed}{\hfill$\Box$\par\medskip\par\relax}
\newcommand{\1}[1]{{\mathbf 1}{\{#1\}}}

\newcommand{\BB}{{\mathcal B}}

\newcommand{\dd}{{\mathtt d}}
\newcommand{\C}{\mathfrak{C}}

\let\phi=\varphi

\newcommand{\lf}{\lfloor}
\newcommand{\rf}{\rfloor}

\newtheorem{theo}{Theorem}[section]
\newtheorem{lm}{Lemma}[section]
\newtheorem{df}{Definition}[section]
\newtheorem{prop}{Proposition}[section]

\allowdisplaybreaks

\title{A conditional quenched CLT for random walks among random 
conductances on $\Z^d$}
\author{Christophe Gallesco$^{~1}$ \and Nina Gantert$^{~2}$ \and 
Serguei Popov$^{~1}$ \and Marina Vachkovskaia$^{~1}$}

\begin{document}

\bibliographystyle{plain}

\maketitle
{\footnotesize

\noindent $^{~1}$  Department of Statistics, 
Institute of Mathematics, Statistics and Scientific Computation,
University of Campinas--UNICAMP, 
rua S\'ergio Buarque de Holanda 651, 13083--859, Campinas SP,
Brazil\\
\noindent e-mail: \texttt{gallesco@ime.usp.br};\\
\noindent e-mail: \texttt{popov@ime.unicamp.br}; 
URL: \texttt{http://www.ime.unicamp.br/$\sim$popov}\\
\noindent e-mail: \texttt{marinav@ime.unicamp.br};
URL: \texttt{http://www.ime.unicamp.br/$\sim$marinav}

\noindent $^{~2}$Fakult\"at f\"ur Mathematik, 
Technische Universit\"at M\"unchen, Boltzmannstr. 3, 85748 Garching, Germany

\noindent e-mail: \texttt{gantert@ma.tum.de}

}

\maketitle

\begin{abstract}
Consider a random walk among random conductances on $\Z^d$ with $d\geq 2$. 
We study the quenched limit law under the usual diffusive scaling
of the random walk conditioned to have its
first coordinate positive. We show that the conditional limit law
is a linear transformation of the product law of a Brownian meander and a $(d-1)$-dimensional
 Brownian motion. 
\\[.3cm] \textbf{Keywords:}  random conductance model, uniform heat kernel bounds,
Brownian meander, reversibility
\\[.3cm] \textbf{AMS 2000 subject classifications:} 60J10, 60K37
\end{abstract}

\section{Introduction and results}
\label{s_intro}
In this paper we study random walks on a $d$-dimensional integer lattice 
with \emph{random conductances}. One can briefly describe the model 
in the following way: initially, weights (i.e., some nonnegative numbers) are 
attached to the edges of the lattice at random. 
The transition probabilities are then defined to be proportional to the weights,
 thus obtaining a reversible Markov chain; due to a well-known 
correspondence between reversible Markov chains and electric networks, 
the weights are also called conductances. We refer to the collection of all 
conductances as ``environment''.
This model attracted considerable attention recently, and, in particular, 
quenched (i.e., for fixed environment) functional central limit theorems and 
heat kernel estimates were obtained in rather general situations, 
see e.g.~\cite{BD,BB,BP,M}
and references therein. We also refer to the survey paper \cite{Bis}.
To prove the quenched functional CLT, one usually uses the so-called 
corrector approach, described in the following way. 
First, one constructs an auxiliary random field (which depends only on 
the environment), with the following property: the sum of 
the corrector and the random walk is a martingale, 
for which it is not difficult to show the CLT. 
Then, using the Ergodic Theorem, one shows that the corrector 
is likely to be small in comparison to the random walk itself.

While this approach has been quite fruitful, it also has its limitations,
 mainly due to the fact that the construction of the corrector 
is not very explicit. For example, it is not clear from this approach
how to prove
the quenched CLT for the random walk with i.i.d.\ conductances in 
\emph{half-space}, even though a similar 
continuous space-and-time problem was solved quite recently~\cite{Rh10}. 
It is therefore important to go beyond the usual setup, proving other types
 of limit laws. 
In this paper, we continue the line of research of~\cite{GP} and~\cite{GP2}
(which were, by their turn, mainly motivated by~\cite{CPSV2,CPSV3}), 
where a one-dimensional model with random conductances 
(but with unbounded jumps) was considered.

We now define the model formally. For $x$, $y\in \Z^d$ with $d\geq 2$, we write 
$x\sim y$ if $x$ and $y$ are neighbors in  the lattice~$\Z^d$ and
we let $\IB_d$ be the set of unordered nearest-neighbor pairs $(x,y)$
of $\Z^d$. Let $(\omega_b)_{b\in \IB_d}$ be non-negative random
variables; $\IP$ stands for the law of this family. We assume that
$\IP$ is stationary and ergodic with respect to the family of
shifts ($\theta_x, x\in \Z^d$). The quantity~$\omega_b$ is usually
called the {\textit{conductance}} of the edge~$b$. The collection
of all conductances $\omega=(\omega_b)_{b\in \IB_d}$ is called the
\textit{environment}. 
If $x\sim y$, we will also write $\omega_{x,y}$ to refer to the conductance 
between~$x$ and~$y$.
For a particular realization $\omega$ of our environment, we define
$\pi_x=\sum_{y\sim x}\omega_{x,y}$. Given that $\pi_x \in (0,\infty)$ 
for all $x\in\Z^d$
(which is $\IP$-a.s.\ the case by Condition~UE below),
the random walk~$X$ in environment~$\omega$ is defined through
its transition probabilities
\[
p_{\omega}(x,y)= \left\{
    \begin{array}{ll}
        \frac{\omega_{x,y}}{\pi_x}, & \mbox{if } y\sim x ,\\
        0, & \mbox{otherwise},
    \end{array}
\right.\
\]
that is, if~$\Po^x$ is the quenched law of the random walk starting
from $x$, we have 
\[
 \Po^x[X(0)=x]=1, \quad \Po^x[X(k+1)=z\mid X(k)=y]=p_{\omega}(y,z).
\]
Clearly, this random walk is $\IP$-a.s.\ reversible with the reversible
measure~$(\pi_x,x\in\Z^d)$.
Also, we denote
by~$\Eo^x$ the quenched expectation for the process starting
from~$x$. When the random walk starts from~$0$, we use the shorter
notations $\Po,\Eo$.

In order to prove our results, we need to make the 
\textit{uniform ellipticity} assumption on the environment:

\textbf{Condition~UE}. There exists $\kappa>0$ such that, 
$\IP$-a.s., $\kappa<\omega_{0,x}<\kappa^{-1}$ for $x\sim 0$.
\medskip

For  all $n\geq 1$, we define the continuous map $(Z^n(t), t\in 
[0,1])$ as the natural polygonal interpolation of the map $k/n\mapsto
n^{-1/2}X(k)$. In other words
\[
 \sqrt{n}Z^n(t) = X(\lfloor nt\rfloor)+(nt-\lfloor
nt\rfloor)X(\lfloor nt\rfloor+1)
\]
with $\lfloor \cdot\rfloor$ the integer part. Also, we denote by~$W^{(d)}=(W_1, \ldots ,W_d)$
the $d$-dimensional standard Brownian motion. Now, let us embed the graph $\Z^d$ in $\R^d$. Denote by $\mathcal{B}=\{{\bf{e}}_1, \dots, {\bf{e}}_d\}$ the canonical basis of $\R^d$ and by $x_1, \dots,x_d$ the vector coordinates in $\R^d$. By Condition UE and as our environment is stationary and ergodic there exists an invertible linear transformation $D:\R^d\to \R^d$ letting the hyperplane $\{x_1=0\}$ invariant and such that the sequence $(DZ^n)_{n\geq 1}$ tends weakly to $W^{(d)}$. Indeed, by Condition UE and ergodicity of the environment, it is well known (cf. \cite{Bis}) that $(Z^n)_{n\geq 1}$ tends weakly to a d-dimensional Brownian motion with a positive definite covariance matrix $\Sigma$. This implies that $\Sigma$ has positive eigenvalues $\lambda_i$ and is diagonalizable in an orthonormal basis. If the law of the environment is also invariant under the symmetries of $\Z^d$, it is known that $\Sigma=\sig^{-1}I$ for some constant $\sig$, where $I$ is the identity matrix.
Thus, there exists a rotation $T$ such that $(TZ^n)_{n\geq 1}$ tends weakly to Brownian motion with diagonal covariance matrix $\Sigma'=(\lambda_i)_{1\leq i\leq d}$ in the basis $\mathcal{B}$. This implies that $((\Sigma')^{-1}TZ^n)_{n\geq 1}$ tends weakly to $W^{(d)}$. Finally, by some unitary transformation $R$, we can rotate the hyperplane $(\Sigma')^{-1}T\{x_1=0\}$ to make it coincide with the hyperplane $\{x_1=0\}$. Now, using the isotropy of $W^{(d)}$ we obtain that $(R(\Sigma')^{-1}TZ^n)_{n\geq 1}$ tends weakly to $W^{(d)}$. For convenience, in the rest of the paper, we will choose $R$ such that $D{\bf{e}}_1 \cdot {\bf{e}}_1>0$.
(R can also involve a reflection). In the case that the law of the environment is also invariant under the symmetries of $\Z^d$, then the last statement is true with $D=\sig^{-1}I$ (where $\sig$ is from the quenched CLT). 

Denoting $X=(X_1,\dots,X_d)$ in the basis $\mathcal{B}$, we define
\[
\hat{\tau}=\inf\{k\geq 1: X_1(k)=0\}
\] 
and 
\[
 \Lambda_n=\{\hat{\tau}>n\}=\{X_1(k)>0\text{ for all }k=1,\ldots,n\}.
\]
Consider the conditional quenched probability measure
$Q_{\omega}^n[\;\cdot\;]:=\Po[~\cdot\mid \Lambda_n]$, for all
$n\geq 1$.  Denote by $C([0,1])$ the space of continuous functions
from $[0,1]$ into $\R^d$. For each~$n$, the random map $DZ^n$ induces a probability measure
$\mu_{\omega}^n$ on $(C[0,1], \BB_1)$:
for
any $A\in \BB_1$,
\[
\mu_{\omega}^n(A):=Q_{\omega}^n[DZ^n\in A].
\]
Let us next recall the formal definition of the
 Brownian meander $W^+$. For this, define $\tau_1=\sup\{s\in [0,1]:
W_1(s)=0\}$ and
$\Delta_1=1-\tau_1$. Then,
\[
W^+(s):=\Delta_1^{-1/2}|W_1(\tau_1+s\Delta_1)|,\phantom{***}0\leq s\leq
1.
\]
We denote by $P_{W^+}\otimes P_{W^{(d-1)}}$ the product law of 
Brownian meander and $(d-1)$-dimensional standard Brownian motion on the time interval $[0,1]$.
Now, we are ready to formulate the quenched invariance principle for
the random walk conditioned to stay positive, which is the main
result of this paper:
\begin{theo}
\label{Theocond}
Under Condition~UE, we have that, $\IP$-a.s.,
 $\mu_{\omega}^n$ tends weakly to $P_{W^+}\otimes P_{W^{(d-1)}}$ as 
$n\to \infty$ (as probability measures on $C[0,1]$).  
\end{theo}

The next result, referred as Uniform Central Limit Theorem (UCLT), will be useful in order to prove Theorem~\ref{Theocond}. Let $W_{\Sigma}$ be a $d$-dimensonal Brownian motion with covariance matrix $\Sigma$ defined above.
Denoting by $\C_b(C([0,1]),\R)$ (respectively,
$\C^u_b(C([0,1]),\R)$) the space of bounded
continuous (respectively, bounded uniformly continuous)
functionals from $C([0,1])$ into $\R$ and by $\mathcal{B}_1$ the Borel
$\sig$-field on $C([0,1])$, we have the following result:
\begin{theo}
\label{Theouni}
Under Condition~UE, the following statements hold 
and are equivalent:
\begin{itemize}
\item[(i)] we have $\IP$-a.s., for all $H>0$ and any 
$F\in \C_b(C([0,1]),\R)$,
\[
\lim_{n \to \infty} \sup_{x\in [-H\sqrt{n},H\sqrt{n}]^d}\Big| 
 {\mathtt E}_{\theta_x \omega}[F(Z^n)]-E[F(W_{\Sigma})]\Big|=0;
\]
\item[(ii)] we have $\IP$-a.s., for all $H>0$ and any
$F\in \C^u_b(C([0,1]),\R)$, 
\[
\lim_{n \to \infty} \sup_{x\in [-H\sqrt{n},H\sqrt{n}]^d}\Big| 
 {\mathtt E}_{\theta_x \omega}[F(Z^n)]-E[F(W_{\Sigma})]\Big|=0;
 \]
 \item[(iii)] we have $\IP$-a.s., for all $H>0$ and any closed 
set~$B$,
\[
\limsup_{n \to \infty} \sup_{x\in [-H\sqrt{n},H\sqrt{n}]^d}
 {\mathtt P}_{\theta_x \omega}[Z^n\in B]\leq P[W_{\Sigma}\in B];
 \]
\item[(iv)] we have $\IP$-a.s., for all $H>0$ and any open set~$G$, 
\[
\liminf_{n \to \infty} \inf_{x\in [-H\sqrt{n},H\sqrt{n}]^d}
{\mathtt P}_{\theta_x \omega}[Z^n\in G]\geq P[W_{\Sigma}\in G];
 \]
 \item[(v)] we have $\IP$-a.s., for all $H>0$ and any $A\in
\mathcal{B}$ such that $P[W_{\Sigma}\in \partial A]=0$,
\[
\lim_{n \to \infty} \sup_{x\in [-H\sqrt{n},H\sqrt{n}]^d}\Big| 
 {\mathtt P}_{\theta_x \omega}[Z^n\in A]-P[W_{\Sigma}\in A]\Big|=0.
 \]
\end{itemize}
\end{theo}

In the next section, we prove some auxiliary results which 
are necessary for the proof of Theorem~\ref{Theocond}. In Section \ref{Sectheouni}, 
we give the proof of Theorem~\ref{Theouni}. Finally, in
Section~\ref{s_proof_Theocond}, we give the proof of Theorem~\ref{Theocond}.
\medskip

We will denote by $C_1$, $C_2$,~$\dots$ the ``global'' 
constants, that is, those that are used all along the paper and by $\gam$, $\gam_1$, $\gam_2$,~$\dots$ the ``local" constants, that is, those
that are used only in the subsection in which they appear for the
first time. For the local constants, we restart the numeration in the
beginning of each subsection. 

Also,
whenever the context is clear, to avoid heavy notations,
 we will not put the integer part symbol $\lf\cdot \rf$. 
For example, for $\delta\in (0,1)$ we will write $X(\delta n)$ 
instead of $X(\lf \delta n\rf)$.

\section{Auxiliary results}
\label{s_aux_results}
In this section, we will prove some technical results that 
will be needed later to prove Theorem~\ref{Theocond}. Instead of considering the process $X$ in the canonical basis $\mathcal{B}$ of $\R^d$ it is also convenient to introduce the embedded graph $\tilde{\Z}^d:=D\Z^d$ with the basis $\mathcal{B'}=\{{\bf{e'}}_1, \dots, {\bf{e'}}_d\}:=D\mathcal{B}$ and consider the process $DX$ in this new basis. All the results obtained in this section concern the original random walk $X$ expressed in $\mathcal {B}$ but they remain valid for $DX$ expressed in $\mathcal{B'}$ with the $\|\cdot\|_1$-norm replaced by the graph distance in $\tilde{\Z}^d$.

Let us introduce the following notations. First, for $a$, $b\in \Z$, $a<b$, 
we denote by $[\![a,b]\!]$ the set $[a,b]\cap \Z$. 
Vectors of $\Z^d$ will be denoted by $x$, $y$ or $z$. For $x\in \Z^d$ 
we denote by $x_1,\dots, x_d$ its coordinates in $\mathcal{B}$.
For $l\in \R$, we denote 
\[
\{l\}_j= \left\{
    \begin{array}{ll}
       \{x=(x_1,\dots, x_d)\in \Z^d: x_j=\lf l \rf\}, & \mbox{if
$l\geq 0$},\\
      \{x=(x_1,\dots, x_d)\in \Z^d: x_j=-\lf l \rf\}, & \mbox{if
$l< 0$},
    \end{array}
\right.\
\]
for $j\in [\![1,d]\!]$. If $F\subset \Z^d$, let us define
\[
\tau_F=\inf\{n\geq 0: X(n)\in F\}\phantom{**} \mbox{and }\phantom{**} \tau^+_F=\inf\{n\geq 1: X(n)\in F\}.
\]
At this point we mention that under
Condition~UE,  we can apply Theorem~1.7 of~\cite{Del} 
to the random walks $Y(n):=X(2n)$ and $Y'(n):=X(2n+1)$, to obtain that 
uniform heat kernel lower and upper bounds are available for this model. That is, 
there exist absolute constants $C_1$, $C_2$, $C_3$ and $C_4$ such that $\IP$-a.s., 
for $n\in \N$,
\begin{equation}
\label{heat_kernel}
p^n_{\omega}(x,y)\leq
\frac{C_1}{n^{d/2}}\exp\Big\{-C_2\frac{\|x-y\|_1^2}{n}\Big\}
\end{equation}
and if $\|x-y\|_1\leq n$ (with $\|\cdot \|_1$ the 1-norm on $\Z^d$) and has the same parity as $n$,
\begin{equation}
\label{heat_kernel2}
p^n_{\omega}(x,y)\geq
\frac{C_3}{n^{d/2}}\exp\Big\{-C_4\frac{\|x-y\|_1^2}{n}\Big\}.
\end{equation}


We denote by $d_1$ the distance induced by the $1$-norm. The heat kernel upper bound (\ref{heat_kernel}) has two simple consequences  gathered in the following
\begin{lm}
\label{Lem0}
Estimate (\ref{heat_kernel}) implies that there exist positive constants $C_5$ and $C_6$ such that $\IP$-a.s., for $h>0$ and $\delta>0$, the following holds.
\begin{itemize}
 \item[(i)] Let $H_1$ and $H_2$ be two parallel hyperplanes in $\Z^d$ orthogonal to ${\bf{e}}_i$ for some $i\in [\![1,d]\!]$ and let us denote by $\mathcal{S}$ the strip delimited by $H_1$ and $H_2$. If $2\leq d_1(H_1, H_2)\leq hn^{1/2}$ then there exists $n_0=n_0(\delta,h)$ such that 
 \[
\sup_{x\in  \mathcal{S}}\Po^x[\tau_{H_1\cup H_2}>\delta^2 n]\leq  C_5\frac{h}{\delta}
 \]
 for all $n\geq n_0$;
 \item[(ii)] Let $x\in \Z^d$. If $A\subset \Z^d$ is such that $d_1(x,A)> hn^{1/2}\geq 1$ then there exists $n_1=n_1(\delta, h)$ such that
\[
\Po^x[\tau_A\leq \delta^2 n]\leq C_6\frac{\delta}{h}
 \]
 for all $n\geq n_1$.
\end{itemize}
\end{lm}
\textit{Proof.} Let us denote by $\mathcal{S}$ the strip delimited by $H_1$ and $H_2$. To prove (i), we just notice that $\Po^x[\tau_{H_1\cup H_2}>\delta^2 n]\leq  \Po^x[X(\delta^2 n)\in \mathcal{S}]$ and apply (\ref{heat_kernel}). More precisely, suppose that $H_1$ and $H_2$ are orthogonal to ${\bf{e}}_1$. With a slight abuse of notation, we also denote by $H_1$ and $H_2$ the coordinates where the hyperplanes $H_1$ and $H_2$ cross the first axis. We have
\begin{align}
\label{ULIU}
\Po^x[X(\delta^2 n)\in \mathcal{S}]
&\leq \sum_{y\in \mathcal{S}}\frac{C_1}{\lf\delta^2n\rf^{d/2}}\exp\Big\{-C_2\frac{\|x-y\|_1^2}{\lf\delta^2n\rf}\Big\}\nonumber\\
&\leq \frac{C_1}{\lf\delta^2n\rf^{d/2}}\sum_{y_1\in [H_1,H_2]}\exp\Big\{-C_2\frac{(y_1-x_1)^2}{\lf\delta^2n\rf}\Big\}\prod_{i=2}^{d}\sum_{y_i\in \Z}\exp\Big\{-C_2\frac{(y_i-x_i)^2}{\lf\delta^2n\rf}\Big\}.
\end{align}
Using (\ref{ULIU}), we can see that there exist positive contants $\gam_1$, $\gam_2$ and $n_0=n_0(\delta,h)$ such that 
\begin{equation*}
\Po^x[X(\delta^2 n)\in \mathcal{S}]\leq \gam_1\int_{0}^{\gam_2\frac{h}{\delta}}\exp{\{-C_2t^2\}}dt
\end{equation*}
for all $n\geq n_0$.
We deduce that there exists a constant $\gam_3>0$ such that 
\begin{equation*}
\Po^x[X(\delta^2 n)\in \mathcal{S}]\leq \gam_3\frac{h}{\delta}
\end{equation*}
for all $n\geq n_0$.

 To prove (ii) we use an argument by Barlow (cf. \cite{Bar} Chapter 3). First, if we denote by $B(x, r)$ the $\|\cdot\|_1$-ball of center $x$ and radius $r:=\lfloor h n^{1/2}\rfloor$ we have that 
$$\Po^x[\tau_A\leq \delta^2 n]\leq \Po^x[\tau_{B^c(x,r)}\leq \delta^2n].$$
Then, we have 
\begin{equation}
\label{FTY}
\Po^x[\tau_{B^c(x,r)}\leq \delta^2n] \leq \Po^x\Big[\|X(\delta^2n)-x\|_1>\frac{r}{2}\Big]+\Po^x\Big[\tau_{B^c(x,r)}\leq \delta^2n, \|X(\delta^2n)-x\|_1 \leq \frac{r}{2}\Big].
\end{equation}
Writing $S=\tau_{B^c(x,r)}$, by the Markov property, the second term of the right-hand side of (\ref{FTY}) equals
\begin{equation*}
\Eo^x\Big[{\bf{1}}_{\{S\leq \delta^2 n\}}\Po^{X_S}\Big[\|X(\lf\delta^2 n\rf-S)-x\|_1\leq \frac{r}{2}\Big]\Big]\leq \sup_{y\in \partial B(x,r+1)}\sup_{m\leq \lf\delta^2n\rf}\Po^y\Big[\|X(\lf\delta^2 n\rf-m)-y\|_1>\frac{r}{2}\Big]
\end{equation*}
where $\partial B(x,r):=\{y\in \Z^d: \|y-x\|_1=r\}$.
Combining this last inequality with (\ref{FTY}) we obtain,
\begin{align*}
\Po^x[\tau_{B^c(x,r)}\leq \delta^2n]
&\leq 2\sup_{y\in \Z^d}\sup_{m\leq \lf\delta^2n\rf}\Po^y\Big[\|X(\lf\delta^2 n\rf-m)-y\|_1>\frac{r}{2}\Big]\nonumber\\
&\leq 2\sup_{y\in \Z^d}\sup_{m\leq \lf\delta^2n\rf}\Po^y\Big[\|X(\lf\delta^2 n\rf-m)-y\|_{\infty}>\frac{r}{2d}\Big]
\end{align*}
where $\|\cdot\|_{\infty}$ is the $\infty$-norm on $\Z^d$. Applying (\ref{heat_kernel}) to bound the last term of the above equation from above and performing the same kind of computations as in the proof of (i), we obtain (ii).
\qed


Next, we prove the following lemma, which gives a uniform lower bound for the probability of progressing in direction ${\bf{e}}_1$ before backstepping to the hyperplane
$\{0\}_1$. 
\begin{lm}
\label{Lem1}
Let $v>0$, then there exist a constant $C_7=C_7(v)>0$ such that we have $\IP$-a.s.,
$\inf_{y\in \{l\}_1}\Po^y[\tau_{\{(v+1)l\}_1}<\tau_{\{0\}_1}]\geq C_7$, for all integers $l\geq 1$.
\end{lm}
\textit{Proof.} 
We are going to show that we can choose $v>0$ small enough in such a way that the statement of Lemma \ref{Lem1} is true for this $v$. The generalization to all $v>0$ is then a direct consequence of the elliptic Harnack inequality. 

For the moment, let $v\in (0,\frac{1}{4})$ and fix $l$ such that $vl\geq 1$.
Then, consider $w\in (v, 1]$. 
We start by writing
\begin{align}
\label{YYY}
\Po^y[\tau_{\{(v+1)l\}_1}<\tau_{\{0\}_1}]
&\geq \Po^y[X_1(wl^2)\geq (v+1)l, \tau_{\{0\}_1}>wl^2]\nonumber\\
&\geq \Po^y[X_1(wl^2)\geq (v+1)l]-\Po^y[ \tau_{\{0\}_1}\leq wl^2].
\end{align}
Next, let us define $\nu:=\lf wl^2\rf$ if $\lf wl^2\rf$ is even or $\nu:=\lf wl^2\rf+1$ otherwise. In the same way, we define  $\rho:=\lf vl\rf$ if $\lf vl\rf$ is even or $\lf vl\rf+1$ otherwise. Observe that in any of these cases,
\begin{equation}
\label{even-even}
\Po^y[X_1(wl^2)\geq (v+1)l]\geq \Po^y[X_1(\nu)> l+\rho].
\end{equation}
We will bound the term of the right-hand side of (\ref{even-even}) from below.
For $y \in \Z^d$, we denote by $\mathcal{P}(y)$ the (non-empty) set of vectors $z\in \Z^d$ that satisfy the following conditions:
$z_1-y_1> \rho$, $\|y-z\|_1$ is even and $\|y-z\|_1 \leq \nu$. 
Applying (\ref{heat_kernel2}), we obtain after some computations
\begin{align}
\label{YYY1}
\Po^y[X_1(\nu)> l+\rho]
&\geq \frac{C_3}{\nu^{d/2}}\sum_{u\in \mathcal{P}(y)} \exp\Big\{-C_4\frac{\|u-y\|_1^2}{\nu} \Big\}\nonumber\\
&\geq \gam_1 \int_{1\wedge 4\sqrt{2}v w^{-1/2}}^{1} \int_{0}^{(1-\frac{u_1}{2})}\dots \int_{0}^{(1-\sum_{i=1}^{d-1}\frac{u_i}{2})}\exp\Big\{-\gam_2 \sum_{i=1}^{d}u_i^2\Big\}du_d\dots du_1\nonumber\\
&:= \gam_1 J(v w^{-1/2})
\end{align}
with $\gam_1$ and $\gam_2$ positive constants depending only on $d$.
By (ii) of Lemma (\ref{Lem0}) we obtain $\Po^y[ \tau_{\{0\}_1}\leq wl^2]\leq C_6w^{1/2}$.
Combining this last inequality with (\ref{YYY}), (\ref{even-even}) and (\ref{YYY1}) we obtain
\begin{equation}
\label{YYY4}
\Po^y[\tau_{\{(v+1)l\}_1}<\tau_{\{0\}_1}] \geq \gam_1 J(vw^{-1/2})-C_6w^{1/2}.
\end{equation}
Observe that for fixed $w$, we have $J(v)\to J(0)>0$ as $v\to 0$, since the integrated function is continuous and positive on its domain of integration and the domain of integration of $J(0)$ has Lebesgue measure bounded from below by $2^{-\frac{d(d-1)}{2}}$. Let 
\[
\eta^*=\max\Big\{\eta\in (0,1]: C_6\eta^{1/2} \leq \frac{1}{4}\gam_1 J(0)\Big\}.
\]
Letting $v< \eta^*\wedge (1/4)$, we can choose a sufficiently small~$w$ in such a way that the
second term of the right-hand side of~(\ref{YYY4}) is smaller than~$\frac{1}{4}\gam_1 J(0)$. 
Once we have chosen~$w$, we can choose~$v$ sufficiently small
in such a way that $J(v)>\frac{J(0)}{2}$. We obtain that
\[
\Po^y[\tau_{\{(v+1)l\}_1}<\tau_{\{0\}_1}] \geq \frac{1}{4}\gam_1 J(0)>0.
\]
This shows Lemma~\ref{Lem1}.
\qed

For $\eps\in (0,1]$, we denote $N:=\lfloor \eps\sqrt{n} \rfloor$. We next prove an upper bound for the probability that the hitting time of the hyperplane ${\{N\}_1}$ is larger than $\eps^{1/2}n$, given
$\Lambda_n$. 
\begin{lm}
\label{Lem2}
There exists a function $f=f(\eps)$ with $\lim_{\eps\to 0}\eps^{-2}f(\eps)=0$ 
such that we have $\IP$-a.s.
\[
\limsup_{n \to \infty}
\Po[\tau_{\{N\}_1}>\eps^{1/2}n\mid \Lambda_n]\leq f(\eps).
\] 
\end{lm}
{\it{Proof.}} Let us begin the proof by sketching the main argument. Consider $\alpha\in (0,1)$, we will show that 
\[
\limsup_{n\to \infty}\Po[\tau_{\{N\}_1}>\eps^{1/2}n\mid \Lambda_n]\leq \limsup_{n\to \infty}\Po[\tau_{\{2^{-1}N\}_1}>\alpha \eps^{1/2}n\mid \Lambda_n]+o_1(\eps)
\] 
when $\eps\to 0$.
\begin{figure}[!htb]
\begin{center}
\includegraphics[scale= 0.5]{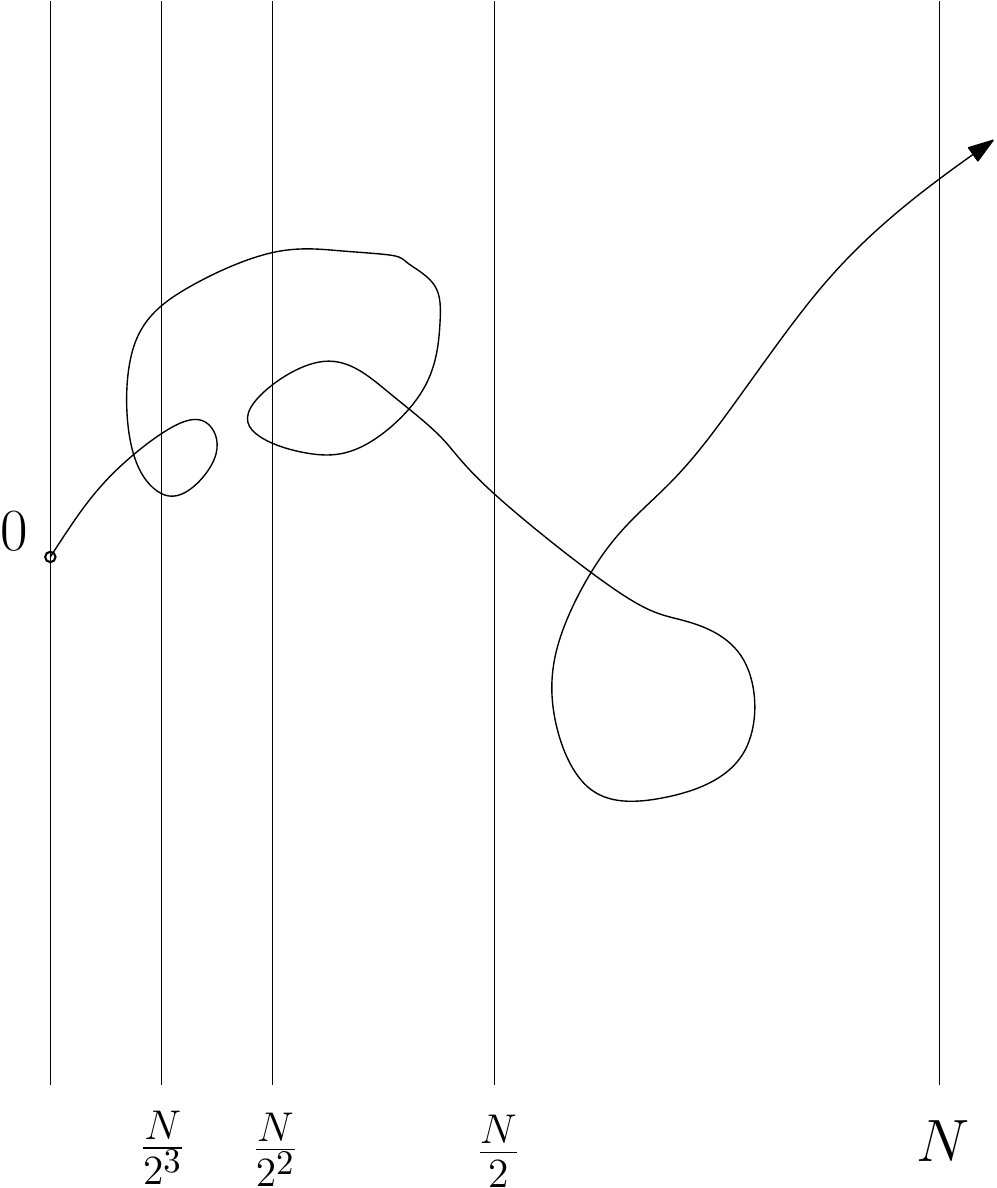}
\caption{Iteration method.}
\label{fig2}
\end{center}
\end{figure}
Then, iterating the argument using hyperplanes of the form $\{2^{-j}N\}_1$ (cf.\ Figure \ref{fig2}) we will have that for all $j\geq 0$,
\[
\limsup_{n\to \infty}\Po[\tau_{\{2^{-j}N\}_1}>\alf^j \eps^{1/2}n\mid \Lambda_n]
\leq \limsup_{n\to \infty}\Po[\tau_{\{2^{-(j+1)}N\}_1}>\alf^{j+1}\eps^{1/2} n
\mid \Lambda_n]+o_j(\eps)
\] 
when $\eps\to 0$.
Finally, restricting $\alf$ to the interval $(\frac{1}{4},1)$, we will show that the $o_j(\eps)$ are 
decreasing fast enough.
Now, let us start the formal argument. Fix $\alf\in (\frac{1}{4},1)$ and let $A_l:=\{\tau_{\{l\}_1}<\tau^+_{\{0\}_1}\}$. We have
\begin{align}
\label{Recu1}
\lefteqn{\Po[\tau_{\{N\}_1}>\eps^{1/2}n\mid \Lambda_n]}
\phantom{***}\nonumber\\
&=\frac{1}{\Po[\Lambda_n]}\Big( \Po[\tau_{\{N\}_1}>
\eps^{1/2}n,\tau_{\{2^{-1}N\}_1}>\alpha \eps^{1/2}
n,
\Lambda_n]+\Po[\tau_{\{N\}_1}>\eps^{1/2}n,\tau_{\{2^{-1}N\}
_1}\leq\alpha\eps^{1/2}n, \Lambda_n]   \Big)\nonumber\\
&\leq \Po[\tau_{\{2^{-1}N\}_1}>\alpha \eps^{1/2}n\mid 
\Lambda_n]+\frac{1}{\Po[\Lambda_n]}\Po[\tau_{\{N\}_1}>\eps^{\frac{1
}{2}}n,\tau_{\{2^{-1}N\}_1}\leq\alpha
\eps^{1/2}n,A_{2^{-1}N}, \Lambda_n].
\end{align}
Then, we have by the Markov property
\begin{align}
\lefteqn{\Po[\tau_{\{N\}_1}>\eps^{1/2}n,\tau_{\{2^{-1}N\}_1}\leq\alpha 
\eps^{1/2}n,A_{2^{-1}N}, \Lambda_n]}\phantom{*******}\nonumber\\
&=\sum_{y\in \{2^{-1}N\}_1}\sum_{k\leq \lfloor \alf \eps^{1/2}n\rfloor}
\Po\Big[X(\tau_{\{2^{-1}N\}_1})=y,\tau_{\{2^{-1}N\}_1}=k,\tau_{\{N\}_1}>\eps^{1/2}n,
A_{2^{-1}N}, \Lambda_n\Big] \nonumber\\
&\leq \max_{y\in \{2^{-1}N\}_1}
\max_{k\leq \lfloor \alf \eps^{1/2}n\rfloor}\Po^y[\tau_{\{N\}_1}
>\eps^{1/2}n-k, \Lambda_{n-k}]\Po[A_{2^{-1}N}].
\end{align}
Now, let us bound from above the term
$\Po^y[\tau_{\{N\}_1}>\eps^{1/2}n-k, \Lambda_{n-k}]$
uniformly in $y\in \{2^{-1}N\}_1$ and in $k\leq \lfloor \alf
\eps^{1/2}n\rfloor$. Observe that, since $\eps\in (0,1]$,
we have
\begin{align}
\Po^y[\tau_{\{N\}_1}>\eps^{1/2}n-k, \Lambda_{n-k}]
&\leq \Po^y[\tau_{\{N\}_1}>(1-\alf)\eps^{1/2}n,\Lambda_{(1-\alf)n}]
\leq \Po^y[\tau_{\{0\}_1\cup\{N\}_1}>(1-\alf)\eps^{1/2}n].
\end{align}

Let $\delta:=\beta^{-1}\eps$, where~$\beta$ is a positive constant
to be determined later. Then, consider~$\eps$ small enough in such
a way that $\delta<(1-\alf)\eps^{1/2}$. Then, divide the time
interval $[0, \lf(1-\alf)\eps^{1/2}n\rf]$ into intervals of size
$\lf \delta^2n \rf$. Denoting 
$S(0,N)=\bigcup_{i=1}^{N-1}\{i\}_1$, we obtain by the Markov
property
\begin{align}
\label{RTY}
\Po^y[\tau_{\{0\}_1\cup\{N\}_1}>(1-\alf)\eps^{1/2}n]
&\leq \Po^y\Big[\tau_{\{0\}_1\cup\{N\}_1}\notin \
\bigcup_{i=1}^{\Big\lf\frac{\lf (1-\alf)\eps^{1/2}n\rf}
 {\lf\delta^{2}n\rf}\Big\rf}((i-1)\lf \delta^2n \rf, i\lf \delta^2n\rf ]\Big]
\nonumber\\
&\leq \Big( \max_{z\in S(0,N)}\Po^z[\tau_{\{0\}_1\cup\{N\}_1}>\delta^2n] 
 \Big)^{(1-\alf)\eps^{1/2}\delta^{-2}-2}
\end{align}
for~large enough $n$.
Using (i) of Lemma \ref{Lem0}, we have for all $z\in S(0,N)$,
\begin{align}
\label{Heat-ker}
\Po^z[\tau_{\{0\}_1\cup\{N\}_1}>\delta^2n] \leq C_5 \frac{\eps}{\delta}
\end{align}
for sufficiently large $n$.
Since $\eps/\delta=\beta$, let us choose the constant $\beta$ such that 
$ C_5 \beta\leq1/2$.
Thus, for~$\eps$ sufficiently small such 
that $\beta^{-1}\eps<(1-\alf)\eps^{1/2}$, we obtain by (\ref{RTY})
\begin{equation*}
\Po^y[\tau_{\{0\}_1\cup\{N\}_1}>(1-\alf)\eps^{1/2}n]\leq 
4\Big(\frac{1}{2}\Big)^{(1-\alf)\eps^{-\frac{3}{2}}\beta^2}.
\end{equation*}
From (\ref{Recu1}), we deduce 
\begin{equation}
\label{Recu2}
\Po[\tau_{\{N\}_1}>\eps^{1/2}n\mid \Lambda_n]\leq 
\Po[\tau_{\{2^{-1}N\}_1}>\alf \eps^{1/2}n\mid \Lambda_n]
+4\Big(\frac{1}{2}\Big)^{ (1-\alf)\eps^{-\frac{3}{2}}\beta^2}
\frac{\Po[A_{2^{-1}N}]}{\Po[\Lambda_n]}.
\end{equation}
Then, we will find an upper bound for the ratio in the second term
of the right-hand side of~(\ref{Recu2}). By the Markov
property we have
\begin{equation}
\label{Recu3}
\frac{\Po[\Lambda_n]}{\Po[A_{2^{-1}N}]}\geq \Po[\Lambda_n\mid A_{2^{-1}N}]
\geq \min_{y\in \{2^{-1}N\}_1}\Po^y[\tau_{\{0\}_1}>n].
\end{equation}
Let $K\geq2\eps$ and let $N'=\lf K\sqrt{n}\rf$.
We start by noting that for any $y\in \{2^{-1}N\}_1$ we have by the
Markov property
\begin{align}
\label{Decomp1}
\Po^y[\tau_{\{0\}_1}>n]
&\geq \Po^y[\tau_{\{0\}_1}>n,\tau_{\{N'\}_1}<\tau_{\{0\}_1}]\geq \min_{z\in\{N'\}_1}\Po^z[\tau_{\{0\}_1}>n]
\Po^y[\tau_{\{N'\}_1}<\tau_{\{0\}_1}].
\end{align}
Let us now bound from below both terms in the right-hand side 
of~(\ref{Decomp1}).

We first show that we can choose a sufficiently large~$K$ in such a
way that $\Po^z[\tau_{\{0\}_1}>n]\geq 1/2$ uniformly in $z\in
\{N'\}_1$. 
Using (ii) of Lemma \ref{Lem0}, we have
$\Po^z[\tau_{\{0\}_1}\leq n]
\leq C_6/K$ for sufficiently large $n$.
Choosing $K$ sufficiently large so that 
$C_6/K \leq1/2$
 we obtain 
\begin{equation}
\label{GLU1}
\Po^z[\tau_{\{0\}_1}>n]\geq \frac{1}{2}
\end{equation}
uniformly in $z\in \{N'\}_1$.
Now going back to equation (\ref{Decomp1}),
 we now show that with probability of order $\eps^{\gam}$ 
with $\gam>0$, starting from the line $\{2^{-1}N\}_1$, 
the random walk reaches the line $\{N'\}_1$  before reaching the line $\{0\}_1$.
By Lemma~\ref{Lem1}, there exists  $C_7>0$ such that 
for every $l>1$, $\Po^u[\tau_{\{2l\}}<\tau_{\{0\}}]\geq C_7$, 
with $u\in \{l\}_1$. Now consider, the following sequence $(U_j)_{j\geq 1}$ 
of hyperplanes defined by 
\[
\left\{
\begin{array}{ll}
U_1&=\{2\lf 2^{-1}N\rf \}_1\\
U_{j+1}&=\{2U_j\}_1.
\end{array}
\right.\
\]
Let $j^{*}$ the smallest~$j$ such that $U_j\geq K\sqrt{n}$.
Using the induction relation, we obtain that for some constant $\gam_1>0$, $j^*\leq
\gam_1 \ln\frac{K}{\eps}$ for large enough $n$. By
convention, set $U_{0}=\{2^{-1}N\}_1$. By the Markov property, we
obtain that uniformly in $y\in \{2^{-1}N\}_1$,
\begin{align}
\label{GLU2}
\Po^y[\tau_{\{N'\}_1}<\tau_{\{0\}_1}]
&\geq \Po^y\Big[\bigcap_{i=1}^{j^*}\{\tau_{U_i}<\tau_{\{0\}_1}\}\Big]
\geq \prod_{i=1}^{j^*}\Big(\min_{u\in U_{i-1}}
\Po^{u}[\tau_{U_i}<\tau_{\{0\}_1}]\Big)  \geq \Big( \frac{\eps}{K}\Big)^{\gam_2}
\end{align}
for some constant $\gam_2>0$ and large enough $n$.
Combining (\ref{Decomp1}), (\ref{GLU1}), and~(\ref{GLU2}) we deduce
\begin{equation}
\label{Recu4}
\min_{y\in \{2^{-1}N\}_1}\Po^y[\tau_{\{0\}_1}>n]
\geq \frac{1}{2}\Big( \frac{\eps}{K}\Big)^{\gam_2}
\end{equation}
for large enough $n$.
Then by (\ref{Recu2}), (\ref{Recu3}) and (\ref{Recu4}) we obtain
\begin{align}
\label{Recu5}
\Po[\tau_{\{N\}_1}>\eps^{1/2}n\mid \Lambda_n]
\leq \Po[\tau_{\{2^{-1}N\}_1}>\alf\eps^{1/2} n
\mid \Lambda_n]+16K^{\gam_2}\eps^{-\gam_2}
\Big(\frac{1}{2}\Big)^{(1-\alf)\eps^{-\frac{3}{2}}\beta^2}.
\end{align}
By the same argument, we can deduce that for all $j\geq 1$ we have 
\begin{equation}
\label{Recu78}
\Po[\tau_{\{2^{-j}N\}_1}>\alf^j \eps^{1/2}n\mid \Lambda_n]
\leq \Po[\tau_{\{2^{-(j+1)}N\}_1}>\alf^{j+1}\eps^{1/2} n
\mid \Lambda_n]+16K^{\gam_2}\Big(\frac{\eps}{2^j}\Big)^{-\gam_2}
\Big(\frac{1}{2}\Big)^{(1-\alf)\beta^2\eps^{-\frac{3}{2}}(4\alf)^j}
\end{equation}
for large enough $n$.
Iterating (\ref{Recu5}) using (\ref{Recu78}), we deduce 
\begin{equation}
\limsup_{n \to \infty}\Po[\tau_{\{N\}_1}>\eps^{1/2}n\mid \Lambda_n]\leq 16K^{\gam_2}
\sum_{j=0}^{\infty} 
\Big(\frac{\eps}{2^j}\Big)^{-\gam_2} 
\Big(\frac{1}{2}\Big)^{(1-\alf)\beta^2\eps^{-\frac{3}{2}}(4\alf)^j}.
\end{equation}
As $\alf \in(\frac{1}{4},1)$, the last series is convergent. 
Define the function $f$ in the statement of Lemma \ref{Lem2} as
\[
f(\eps):=16K^{\gam_2} \sum_{j=0}^{\infty} \Big(\frac{\eps}{2^j}\Big)^{-\gam_2}
 \Big(\frac{1}{2}\Big)^{(1-\alf)\beta^2\eps^{-\frac{3}{2}}(4\alf)^j}.
\]
Using the dominated convergence theorem, it is straightforward to
 show that $\lim_{\eps\to 0} \eps^{-2} f(\eps)=0$. This proves
Lemma~\ref{Lem2}.
\qed

In the next lemma, $N$ still stands for $\lf \eps \sqrt{n} \rf$. However, the quantities 
(like $\alpha$, $\delta$,  $\beta$, ...)  defined in the proof of the lemma are not related
to the corresponding quantities defined in the proof of Lemma \ref{Lem2}. The next lemma controls the ``transversal fluctuations'' of $X_2, \ldots ,X_d$, given
$\Lambda_n$.

\begin{figure}[!htb]
\begin{center}
\includegraphics[scale= 0.5]{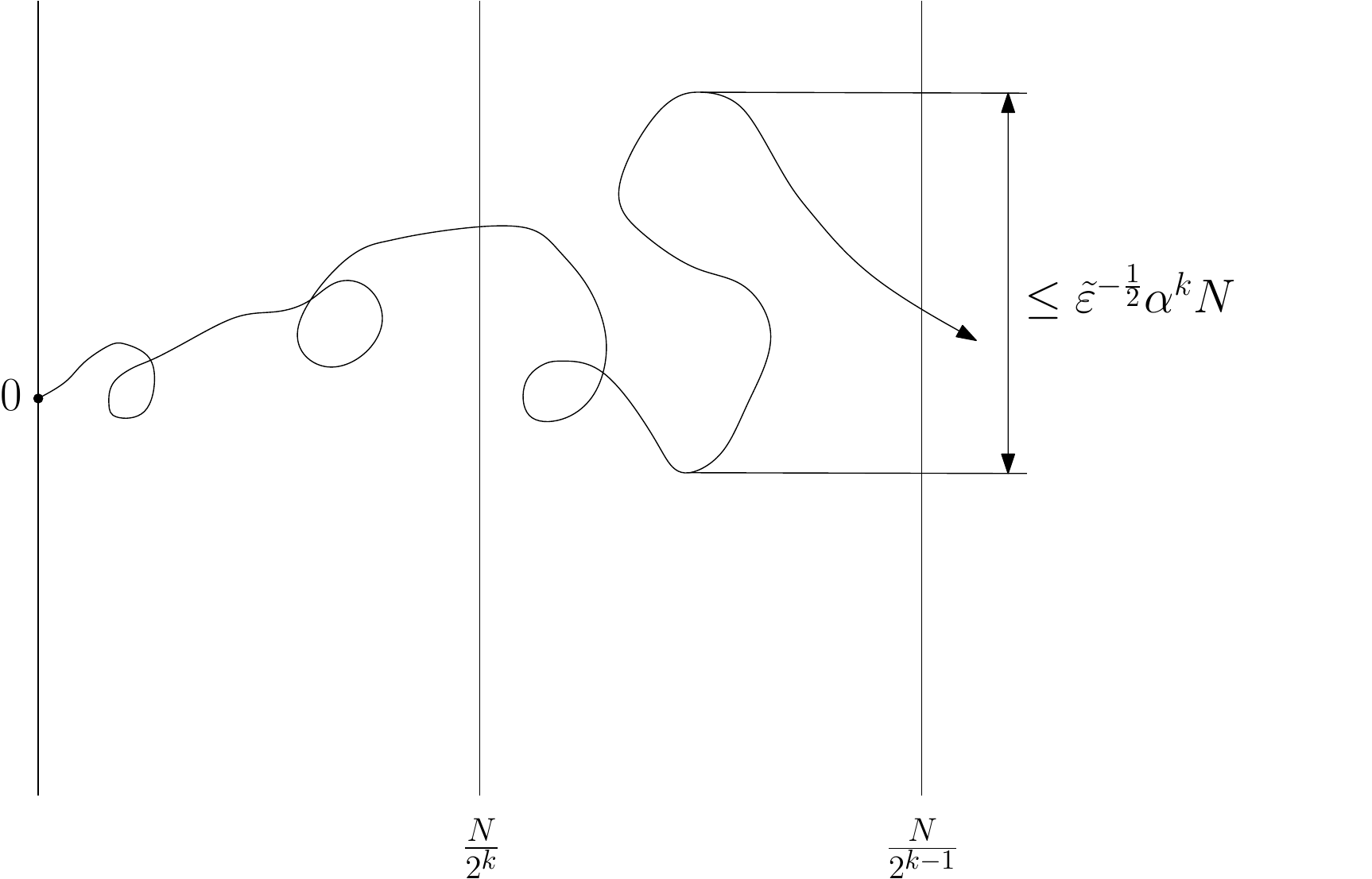}
\caption{On the definition of $G_k$.}
\label{fig3}
\end{center}
\end{figure}

\begin{lm}
\label{Lem3}
We have $\IP$-a.s.,
\[
\limsup_{n\to \infty}\Po\Big[\max_{i\in [\![2,d]\!]}
\sup_{j\leq \tau_{\{N\}_1}}|X_i(j)|
>\eps^{-1/2}N\mid\Lambda_n\Big]\leq g(\eps)
\]
with $\lim_{\eps \to 0}\eps^{-2}g(\eps)=0$.
\end{lm}
{\it Proof.} First, observe that, by symmetry, it suffices to show
that there exists $g'=g'(\eps)$ such that
\begin{equation}
\label{Lem32}
\limsup_{n\to \infty}\Po\Big[\sup_{j\leq \tau_{\{N\}_1}}|X_i(j)|
>\eps^{-1/2}N\mid\Lambda_n\Big]\leq g'(\eps)
\end{equation}
with $\lim_{\eps \to 0}\eps^{-2}g'(\eps)=0$ for some
 $i\in [\![2,d]\!]$. For the sake of simplicity, let us take $i=2$
in the rest of the proof.
Fix $\alf \in (\frac12,1)$ and let $\tilde{\eps}^{-1/2}:=\frac{1-\alf}{\alf}\eps^{-1/2}>2$. We introduce the
following sequence of events (cf. Figure \ref{fig3}),
\[
G_k=\Big\{ \sup_{j\in(\tau_{\{2^{-k}N\}_1}, 
\tau_{\{2^{-k+1}N\}_1}]}|X_2(j)-X_2(\tau_{\{2^{-k}N\}_1})|
\leq \tilde{\eps}^{-1/2} \alf^k N \Big\}
\]
for $k\geq 1$, with the convention that $\sup_{j\in\emptyset}\{\cdot\}=0$.
Then, we denote
\[B_k^{\delta}=\{\tau_{\{2^{-k}N\}_1}\leq \delta n\}\cap \{\tau_{\{2^{-k}N\}_1}<\tau_{\{0\}_1}\}\]
for $\delta\in (0,1]$ and $k\geq 1$.

Now, observe that on the event 
$B_0^{\delta}\cap (\cap_{k\geq 1}G_k)$ we have that $\sup_{j\leq
\tau_{\{N\}_1}}|X_2(j)|\leq \eps^{-1/2}N$ since
$\alf \in (\frac12,1)$. This implies that
\begin{equation*}
\Po\Big[\sup_{j\leq \tau_{\{N\}_1}}|X_2(j)|
\leq \eps^{-1/2}N\mid \Lambda_n\Big]
\geq \Po\Big[ B_0^{\delta}\cap (\bigcap_{k\geq 1}G_k) \mid \Lambda_n\Big].
\end{equation*}
In order to prove Lemma \ref{Lem3}, we will show that 
$\liminf_{n\to \infty}\Po[B_0^{\delta}\cap (\cap_{k\geq 1}G_k) \mid \Lambda_n]$ tends to 1 
when $\eps \to0$. We start by writing
\begin{align}
\label{FG41}
\Po\Big[B_0^{\delta}\cap (\bigcap_{k\geq 1}G_k) \mid \Lambda_n\Big]
&=\Po[B_0^{\delta}\mid \Lambda_n]-\Po\Big[B_0^{\delta}
\cap (\bigcap_{k\geq 1}G_k)^c\mid \Lambda_n\Big]\nonumber\\
&\geq \Po[B_0^{\delta}\mid \Lambda_n]-\sum_{k=1}^{\lf \frac{\ln N}{\ln 2}\rf}
\Po[B_0^{\delta}\cap G_k^c\mid \Lambda_n].
\end{align}
From now on, we dedicate ourselves to bounding from above the terms $\Po[B_0^{\delta}\cap G_k^c\mid \Lambda_n]$ for $k\leq \lf \frac{\ln N}{\ln 2}\rf$. We have by the Markov property,
\begin{align*}
\Po[B_0^{\delta}\cap G_k^c\mid \Lambda_n]
&\leq \Po[B_k^{\delta}\cap G_k^c\mid \Lambda_n]=\frac{1}{\Po[\Lambda_n]}\Po[B_k^{\delta},G_k^c, \Lambda_n]\nonumber\\
&=\frac{1}{\Po[\Lambda_n]}\sum_{j\leq \lf \delta n \rf}
\sum_{y\in \{2^{-k}N\}_1}\Po\Big[B_k^{\delta},G_k^c, \Lambda_n, 
\tau_{\{2^{-k}N\}_1}=j, X(\tau_{\{2^{-k}N\}_1})=y\Big]\nonumber\\
&\leq \frac{\Po[B^1_k]}{\Po[\Lambda_n]}\max_{j\leq \lf \delta n \rf}
\max_{y\in \{2^{-k}N\}_1}
\Po^y\Big[\sup_{i\leq \tau_{\{2^{-k+1}N\}_1}}|(X(i)-y)\cdot {\bf{e}}_2|
>\tilde{\eps}^{-1/2}\alf^kN, \Lambda_{n-j}\Big]\nonumber\\
&\leq \frac{\Po[B^1_k]}{\Po[\Lambda_n]}\max_{y\in \{2^{-k}N\}_1}
\Po^y\Big[\sup_{i\leq \tau_{\{2^{-k+1}N\}_1}}|(X(i)-y)\cdot {\bf{e}}_2|
>\tilde{\eps}^{-1/2}\alf^kN, \Lambda_{(1-\delta)n}\Big].
\end{align*}
Using again the Markov property, we obtain
\begin{equation*}
\frac{\Po[\Lambda_n]}{\Po[B^1_k]}\geq \Po[\Lambda_n\mid B^1_k]\geq 
 \min_{y\in \{2^{-k}N\}_1}\Po^y[\tau_{\{0\}_1}>n].
\end{equation*}
By the same argument which we used in Lemma~\ref{Lem2} to treat the term
$\min_{y\in \{2^{-1}N\}_1}\Po^y[\tau_{\{0\}_1}>n]$ 
(cf.\ the derivation of~\eqref{Recu4}), we obtain, for large enough $n$ and all $k\leq \lf \frac{\ln N}{\ln 2}\rf$,
\begin{equation}
\label{JKL1}
\frac{\Po[B^1_k]}{\Po[\Lambda_n]}\leq  \gam_1 \Big(\frac{K2^k}{\eps} \Big)^{\gam_2}
\end{equation}
for some positive constants $\gam_1$, $\gam_2$ and $K$ from Lemma~\ref{Lem2}.
Now, we need to bound the terms 
$$\Po^y\Big[\sup_{i\leq \tau_{\{2^{-k+1}N\}_1}}|(X(i)-y)\cdot
{\bf{e}}_2|>\tilde{\eps}^{-1/2}\alf^kN, \Lambda_{(1-\delta)n}\Big]$$ from
above, uniformly in $y\in \{2^{-k}N\}_1$. In order not to carry on
heavy notations we treat the case $y_2=0$. However, as one can
check, the bound we will obtain is uniform in $y\in \{2^{-k}N\}_1$.
Let 
\[
 E_k=\{(x_1,\dots, x_d)\in \Z^d: x_2=\pm \lf \tilde{\eps}^{-1/2}\alf^kN\rf\}.
\]
 We start by writing
\begin{align}
\label{POICA}
\Po^y\Big[\sup_{i\leq \tau_{\{2^{-k+1}N\}_1}}|X_2(i)|>\tilde{\eps}^{-1/2}\alf^kN,
 \Lambda_{(1-\delta)n}\Big]
&= \Po^y[\tau_{E_k}<\tau_{\{2^{-k+1}N\}_1}, \tau_{\{0\}_1}>(1-\delta)n]
\nonumber\\
&\leq \Po^y[\tau_{E_k}<\tau_{\{2^{-k+1}N\}_1\cup \{0\}_1}]
+\Po^y[\tau_{E_k}>(1-\delta)n].
\end{align}
Let us bound the first term of the right-hand side of~(\ref{POICA}) from above. 
To do so, we first write
\begin{align}
\label{OKP}
 \lefteqn{\Po^y[\tau_{E_k}<\tau_{\{2^{-k+1}N\}_1\cup \{0\}_1}]}
\phantom{*******}\nonumber\\
 &\leq \Po^y[\tau_{\{\tilde{\eps}^{-1/2}\alf^kN\}_2}<\tau_{\{2^{-k+1}N\}_1\cup \{0\}_1}]
+\Po^y[\tau_{\{-\tilde{\eps}^{-1/2}\alf^kN\}_2}<\tau_{\{2^{-k+1}N\}_1\cup \{0\}_1}].
\end{align}
We treat the first term of the right-hand side of~(\ref{OKP}) (the
method for the second term is similar). Let $L\in
(2,\tilde{\eps}^{-1/2})$ and divide the interval $[0, \lf
\tilde{\eps}^{-1/2}\alf^kN\rf]$ into intervals of size $\lf L2^{-k}N\rf$. Furthermore, let
\[
F_k=\bigcup_{j=1}^{\lf 2^{-k+1}N\rf-1}\{j\}_1.
\]
We have by the Markov property,
\begin{align}
\label{DFG}
\lefteqn{\Po^y[\tau_{\{\eps^{-1/2}\gam\alf^kN\}_2}<\tau_{\{2^{-k+1}N\}_1\cup 
\{0\}_1}]}\phantom{**********}\nonumber\\
&\leq \Po^y\Big[\bigcap_{j=1}^{\Big \lf\frac{\lf \tilde{\eps}^{-1/2}\alf^k N\rf}
{\lf L2^{-k}N\rf}\Big \rf}\{\tau_{\{j\lf L2^{-k}N\rf\}_2}<\tau_{\{2^{-k+1}N\}_1
\cup \{0\}_1}\}\Big]\nonumber\\
&\leq \prod_{j=1}^{\lf L^{-1}\tilde{\eps}^{-1/2}(2\alf)^k \rf-2}
\max_{z\in \{(j-1)\lf L2^{-k}N\rf\}_2\cap F_k}\Po^z[\tau_{\{j\lf L2^{-k}N\rf\}_2}
<\tau_{\{2^{-k+1}N\}_1\cup \{0\}_1}].
\end{align}
Let us show that 
\[\max_{z\in \{(j-1)\lf L2^{-k}N\rf\}_2\cap F_k}
\Po^z[\tau_{\{j\lf L2^{-k}N\rf\}_2}<\tau_{\{2^{-k+1}N\}_1\cup \{0\}_1}]
\leq\frac{1}{2}\]
 for $\eps$ sufficiently small and $L$ sufficiently large belonging to 
$(2,\tilde{\eps}^{-1/2})$. Consider $w\in (4,L^2)$, we have for $z\in \{(j-1)\lf L2^{-k}N\rf\}_2\cap F_k$,
\begin{align}
\label{RTYU1}
\Po^z[\tau_{\{j\lf L2^{-k}N\rf\}_2}>\tau_{\{2^{-k+1}N\}_1\cup \{0\}_1}]
&\geq \Po^z[\tau_{\{2^{-k+1}N\}_1\cup \{0\}_1}\leq w2^{-2k}N^2, 
\tau_{\{j\lf L2^{-k}N\rf \}_2}>w2^{-2k}N^2]\nonumber\\
&\geq \Po^z[\tau_{\{2^{-k+1}N\}_1\cup \{0\}_1}\leq w2^{-2k}N^2]-
\Po^z[ \tau_{\{j\lf L2^{-k}N\rf\}_2}\leq w2^{-2k}N^2].
\end{align}
Using (i) of Lemma \ref{Lem0}, we deduce
\begin{align}
\label{RTYU2}
\Po^z[\tau_{\{2^{-k+1}N\}_1\cup \{0\}_1}\leq w2^{-2k}N^2]
\geq 1- C_5w^{-1/2}.
\end{align}
Using (ii) of Lemma \ref{Lem0}, we obtain for all $j\geq 1$,
\begin{align}
\label{RTYU3}
\Po^z[ \tau_{\{j\lf L2^{-k}N\rf\}_2}\leq w2^{-2k}N^2] \leq C_6 \frac{w^{1/2}}{L}.
\end{align}
Combining (\ref{RTYU1}), (\ref{RTYU2}) and (\ref{RTYU3}) we obtain for all $j\geq 1$,
\begin{equation}
\Po^z[\tau_{\{jL2^{-k}N\}_2}>\tau_{\{2^{-k+1}N\}_1\cup \{0\}_1}]
\geq 1- C_5w^{-1/2}- C_6 \frac{w^{1/2}}{L}  .
\end{equation}
First, choose $w$ sufficiently large such that 
$C_5w^{-1/2}\leq 1/4$ 
and thus choose $L$ sufficiently large in such a way 
that $C_6 w^{1/2}/L \leq 1/4$. 
We obtain 
\begin{equation}
\label{ISItron}
\Po^z[\tau_{\{jL2^{-k}N\}_2}>\tau_{\{2^{-k+1}N\}_1\cup \{0\}_1}]\geq \frac{1}{2}.
\end{equation}
Now using (\ref{OKP}), (\ref{DFG}) and (\ref{ISItron}) we have since $\tilde{\eps}^{-1/2}>L$,
\begin{equation}
\label{JKL2}
\Po^y[\tau_{E_k}<\tau_{\{2^{-k+1}N\}_1\cup \{0\}_1}]
\leq 16\Big(\frac{1}{2}\Big)^{\lf L^{-1}\tilde{\eps}^{-1/2}(2\alf)^k \rf}.
\end{equation}
Next, let us treat the term $\Po^y[\tau_{E_k}>(1-\delta)n]$. Let
$\eta=\beta^{-1}\tilde{\eps}$ where $\beta$ is a positive constant to be chosen
later. Then suppose that $\eps$ is sufficiently small such that
$\eta \tilde{\eps}^{-1/2} \alf^k<1-\delta$ and divide the time interval
$[0,\lf(1-\delta)n\rf]$ into intervals of size $\lf
\eta^2\tilde{\eps}^{-1}\alf^{2k} n\rf$. Denoting by 
\[
H(E_k)=\bigcup_{j=-\lf \tilde{\eps}^{-1/2}\alf^{k}N\rf+1}^{\lf
 \tilde{\eps}^{-1/2}\alf^{k}N\rf-1}\{j\}_2,
\]
we obtain by the Markov property
\begin{align}
\label{RETYU}
\Po^y[\tau_{E_k}>(1-\delta)n]
&\leq \Po^y\Big[\tau_{E_k}\notin \bigcup_{i=1}^{\Big\lf\frac{\lf (1-\delta)n\rf}
{\lf\eta^2\tilde{\eps}{-1} \alf^{2k}n\rf}\Big\rf}((i-1)
\lf \eta^2\tilde{\eps}^{-1}\alf^{2k}n\rf, i\lf \eta^2\tilde{\eps}^{-1}\alf^{2k}n\rf]\Big]
\nonumber\\
&\leq \Big(\max_{z\in H(E_k)}\Po^z[\tau_{E_k}>\eta^2\eps^{-1}\alf^{2k}n] 
\Big)^{ (1-\delta)(\eta\tilde{\eps}^{-1/2}\alf^k)^{-2}-2}
\end{align}
for $n$ sufficiently large.
We now bound the term $\Po^z[\tau_{E_k}>\eta^2\tilde{\eps}^{-1}\alf^{2k}n]$ from above
 uniformly in $z\in H(E_k)$. Using (i) of Lemma \ref{Lem0}, we have 
\begin{align*}
\Po^z[\tau_{E_k}>\eta^2\tilde{\eps}^{-1}\alf^{2k}n]\leq C_5 \frac{\tilde{\eps}}{\eta}.
\end{align*}
Since $\tilde{\eps}\eta^{-1}=\beta$, choose $\beta$ small enough such that $C_5\beta\leq 1/2$.
 For $\eps$ sufficiently small such that
  $\eta \tilde{\eps}^{-1/2} \alf^k<1-\delta$, we obtain using~(\ref{RETYU}),
\begin{equation}
\label{JKL3}
\Po^y[\tau_{E_k}>(1-\delta)n]
\leq 4 \Big(\frac{1}{2}\Big)^{(1-\delta)(\beta^{-1}\tilde{\eps}^{1/2}\alf^k)^{-2}}.
\end{equation}
Combining (\ref{POICA}), (\ref{OKP}), (\ref{JKL2}) and (\ref{JKL3}), we deduce that, $\IP$-a.s., for all large enough $n$ and $k\leq \lf \frac{\ln N}{\ln 2}\rf$,
\begin{align}
\label{ULYS1}
\lefteqn{\max_{y\in \{2^{-k}N\}_1}
\Po^y\Big[\sup_{i\leq \tau_{\{2^{-k+1}N\}_1}}|(X(i)-y)\cdot {\bf{e}}_2|
>\tilde{\eps}^{-1/2}\alf^kN, \Lambda_{(1-\delta)n}\Big]}\phantom{****************************}\nonumber\\
&\leq 16 \Big(\frac{1}{2}\Big)^{L^{-1}\tilde{\eps}^{-1/2}(2\alf)^{k}}+4\Big(\frac{1}{2}\Big)^{(1-\delta)(\beta^2\tilde{\eps}^{-1}\alf^{-2k})}.
\end{align}
Using (\ref{JKL1}) and (\ref{ULYS1}), we obtain for all large enough $n$ and $k\leq \lf \frac{\ln N}{\ln 2}\rf$,
\begin{equation*}
\Po[B_0^{\delta}\cap G_k^c\mid \Lambda_n]\leq
 \gam_1K^{\gam_2} 2^{k\gam_2+1}\eps^{-\gam_2}
\Big(16\Big(\frac{1}{2}\Big)^{L^{-1}\tilde{\eps}^{-1/2}(2\alf)^k }
+4\Big(\frac{1}{2}\Big)^{(1-\delta)(\beta^{2}\tilde{\eps}^{-1}\alf^{-2k})}\Big).
\end{equation*}
We finally deduce that, $\IP$-a.s., for large enough $n$,
\begin{equation*}
\sum_{k=1}^{\lf \frac{\ln N}{\ln 2}\rf}
\Po[B_0^{\delta}\cap G_k^c\mid \Lambda_n]\leq
\sum_{k=1}^{\infty} \gam_1K^{\gam_2} 2^{k\gam_2+1}\eps^{-\gam_2}
\Big(16\Big(\frac{1}{2}\Big)^{L^{-1}\tilde{\eps}^{-1/2}(2\alf)^k }
+4\Big(\frac{1}{2}\Big)^{(1-\delta)(\beta^{2}\tilde{\eps}^{-1}\alf^{-2k})}\Big).
\end{equation*}
Observe that since $\alf\in (\frac{1}{2},1)$, the series above converges. Let $\delta=\eps^{1/2}$, we have for $\eps<1/4$,
\begin{equation*}
\sum_{k=1}^{\lf \frac{\ln N}{\ln 2}\rf}
\Po[B_0^{\delta}\cap G_k^c\mid \Lambda_n]\leq 
 \sum_{k=1}^{\infty} \gam_1K^{\gam_2} 2^{k\gam_2+1}\eps^{-\gam_2}
\Big(16\Big(\frac{1}{2}\Big)^{ \frac{1-\alf}{\alf}L^{-1}\eps^{-1/2}(2\alf)^k }
+4\Big(\frac{1}{2}\Big)^{1/2( \frac{1-\alf}{\alf})^2\beta^{2}\eps^{-1}\alf^{-2k})}\Big).
\end{equation*}
Let
\[
h(\eps):= \sum_{k=1}^{\infty} \gam_1K^{\gam_2} 2^{k\gam_2+1}\eps^{-\gam_2}
\Big(16\Big(\frac{1}{2}\Big)^{ \frac{1-\alf}{\alf}L^{-1}\eps^{-1/2}(2\alf)^k }
+4\Big(\frac{1}{2}\Big)^{1/2(\frac{1-\alf}{\alf})^2\beta^{2}\eps^{-1}\alf^{-2k})}\Big).
\]
By the Lebesgue dominated convergence theorem, we have $\eps^{-2}h(\eps)\to 0$ 
as $\eps\to 0$.
Using~(\ref{FG41}) and Lemma~\ref{Lem2} (since $\delta = \eps^{1/2}$) we have for $\eps<1/4$,
\begin{equation*}
\liminf_{n \to \infty}\Po\Big[B_0^{\delta}\cap (\bigcap_{k\geq 1}G_k) 
\mid \Lambda_n\Big]\geq 1-f(\eps)-h(\eps).
\end{equation*}
This last term tends to 1 as $\eps\to 0$. Now, take 
$g'(\eps):=f(\eps)+h(\eps)$ to show (\ref{Lem32}) and therefore Lemma~\ref{Lem3}. 
\qed


\section{Proof of the UCLT}
\label{Sectheouni}

In this section we prove Theorem~\ref{Theouni}. The proof is 
similar in spirit to the proof of Theorem~1.2 of~\cite{GP},
nevertheless it is greatly simplified in the present case by the
use of the heat kernel upper bounds. In order to take advantage of the natural left shift on the space $C(\R_+)$ of continuous functions from $\R_+$ into $\R^d$, we will rather prove Theorem \ref{Theouni} for $Z^n$ assuming values in $C(\R_+)$ instead of $C([0,1])$. Then, the result for $Z^n$ assuming values in $C([0,1])$ will be easily obtained by the mapping theorem (cf.\ \cite{Billing}). Let $\C^u_b(C(\R_+),\R)$ be
the space of bounded
uniformly continuous functionals from $C(\R_+)$ into~$\R$. In this section, we write ~$W$
for the $d$-dimensional Brownian motion with covariance matrix $\Sigma$ from section \ref{s_intro}. The first step is to prove the
following 
\begin{prop}
\label{UCLTprop1}
For all $F\in \C^u_b(C(\R_+),\R)$, we have $\IP$-a.s., for every 
$H>0$, 
\begin{equation*}
\lim_{n\to \infty}\sup_{x\in [-H\sqrt{n},H\sqrt{n}]^d}\Big| 
{\mathtt E}_{\theta_x \omega}[F(Z^{n})]-E[F(W)]\Big|=0.
\end{equation*}
\end{prop}
 \medskip
 
 \noindent
Fix $F\in \C^u_b(C(\R_+),\R)$. We will prove that, $\IP$-a.s., 
for every $\tilde{\eps}, H>0$, 
\begin{equation}
\label{UCLTrefor}
\sup_{x\in [-H\sqrt{n},H\sqrt{n}]^d}\Big| 
{\mathtt E}_{\theta_x \omega}[F(Z^{n})]-E[F(W)]\Big|\leq
\tilde{\eps}
\end{equation}
for~$n$ large enough. Before this, we need to introduce some definitions 
and prove an intermediate result. 
Let~$\dd$ be the distance on the space $C_{\R_+}$ 
defined by
\[
\dd(f,g)=\sum_{n=1}^{\infty}2^{-n+1}
\min\Big\{1,\sup_{s\in [0,n]}\|f(s)-g(s)\|\Big\}
\]
with $\| \cdot \|$ the euclidian norm on $\R^d$.
Now, for any given $\eps>0$, let 
\begin{equation}
\label{UUFMW}
h_{\eps}:=\max\Big\{h\in (0,1]: P\Big[\sup_{s\leq h}\|W(s)\|>\eps\Big]
+P\Big[\sup_{s\leq h}\dd(\theta_sW,W)>\eps\Big]\leq
\frac{\eps}{2}\Big\}.
\end{equation}
Observe that $h_{\eps}>0$ for $\eps>0$ and $h_{\eps}\to 0$ when $\eps\to 0$.
Next, adapting section 3 of~\cite{GP} we introduce the following
\begin{df}
For a given realization of the environment $\omega$ and $N\in \N$, we say that $x\in \Z^d$ is $(\eps,
 N)$-good,
 if 
\begin{itemize}
\item
$\min \Big \{n\geq 1: \big| 
{\mathtt E}_{\omega}[F(Z^{m})]-E[F(W)]\big| \leq \eps,
\phantom{*}\mbox{for all $m\geq n$}\Big \}\leq N$;
\item
$ \mathtt{P}_{\theta_x \omega}\Big[\sup_{s\leq h_{\eps}}\|Z^{m}(s)\|
\leq \eps,\sup_{s\leq
h_{\eps}}\dd(\theta_sZ^{m},Z^{m})\leq \eps\Big] 
\geq 1-\eps$, for
all $m\geq N$.
\end{itemize}
\end{df}
We now show that starting from a site $x\in [-H\sqrt{n},H\sqrt{n}]^d$, 
with high probability, the random walk $X$ will meet a
$(\eps,n)$-good site 
 at a distance at most
 $h'\sqrt{n}$ before time~$hn$ (unlike as in ~\cite{GP}, there is no need here 
to introduce the notion of a \textit{nice} site since by~(\ref{heat_kernel}), 
every point in $[-H\sqrt{n}, H\sqrt{n}]^d$ is nice).
We denote by 
$\mathcal{G}$ the set of
$(\eps,n)$-good sites in $\Z^d$.
\begin{prop}
\label{reacgoodsite}
 Fix $h'>0$. For any $\eps_1>0$, we can 
choose $\eps$ small enough
in such a way that we have $\IP$-a.s., for all sufficiently large~$n$ and all $x\in [-H\sqrt{n},H\sqrt{n}]^d$:
\begin{itemize}
\item[(i)] 
$\Po^x[\tau_{\mathcal{G}}> h_{\eps}n]\leq \eps_1;$
\item[(ii)] 
$\Po^x\Big[ \sup_{j\leq h_{\eps}n}\|X(j)-X(0)\|>
h'\sqrt{n}\Big]\leq \eps_1.$
\end{itemize}
\end{prop}

\noindent\textit{Proof.}
Fix  $\eps$. Then, for any $\eps'>0$ 
there exists~$N$ such that
\[
 \IP[0\text{ is $(\eps,N)$-good}] 
      > 1-\eps'.
\]
By the Ergodic
Theorem, we have $\IP$-a.s.\ for all $n>n_1(\omega)$,
\begin{equation}
\big|\{x \in [-2H\sqrt{n}, 2H\sqrt{n}]^d 
\text{ and $x$ is not
$(\eps,N)$-good}\}\big| <
5^d\eps'H^dn^{\frac{d}{2}}.
\end{equation}
Let us define 
$$\mathsf{Bad}:=\{x \in [-2H\sqrt{n}, 2H\sqrt{n}]^d 
\text{ and $x$ is not
$(\eps,N)$-good}\}$$
and $\mathsf{Cub}:=  [-2H\sqrt{n},2H\sqrt{n}]^d$.

In order to show (i) we observe that for all $x\in [-H\sqrt{n}, H\sqrt{n}]^d$,
\begin{align}
\label{OKJ}
\Po^x[\tau_{\mathcal{G}}>  h_{\eps}n]
\leq \Po^x[X(h_{\eps}n)\in \mathsf{Bad}]+\Po^x[\tau_{ \mathsf{Cub}^c}\leq h_{\eps}n].
\end{align}
For the second term of the right-hand side of (\ref{OKJ}), we apply 
(ii) of Lemma \ref{Lem0} to obtain that $\Po^x[\tau_{ \mathsf{Cub}^c}\leq h_{\eps}n]\leq \gam_2 h_{\eps}$. Thus, we can choose $\eps$ small enough in such a way that $\Po^x[\tau_{ \mathsf{Cub}^c}\leq h_{\eps}n]\leq \eps_1/2$.
Then, using (\ref{heat_kernel}) and the fact that $|\mathsf{Bad}|<5^d\eps'H^dn^{\frac{d}{2}}$ for large $n$, we can show that uniformly 
in $x\in \mathsf{Bad}\cap [-H\sqrt{n},H\sqrt{n}]^d$ we have 
$\Po^x[X(h_{\eps}n)\in \mathsf{Bad}]\leq \gam_1 \eps'/h_{\eps}$ for $n$ sufficiently large.
Thus, choosing $\eps'$ sufficiently small in such a way 
that $\gam_1 \eps'/h_{\eps}\leq \eps_1/2$ we obtain $\Po^x[X(h_{\eps}n)\in \mathsf{Bad}]\leq\eps_1/2$. 

To show (ii), we notice that 
\begin{equation}
\label{UHG}
\Po^x\Big[ \sup_{j\leq h_{\eps}n}\|X(j)-X(0)\|>
h'\sqrt{n}\Big] =\Po^x[\tau_{B^c(x,h'\sqrt{n})}\leq h_{\eps}n]
\end{equation}
with $B(x,r)$ the euclidian ball of center $x$ and radius $r$. Now, we can apply (ii) of Lemma \ref{Lem0} to the right-hand term of (\ref{UHG}) to obtain that 
\begin{equation*}
\Po^x\Big[
\sup_{j\leq h_{\eps}n}\|X(j)-X(0)\|>
h'\sqrt{n}\Big]\leq \gam_3\frac{h_{\eps}^{1/2}}{h'}.
\end{equation*}
Finally, choosing $\eps$ sufficiently small such that $\gam_3 h^{1/2}/h'\leq \eps_1$ we obtain (ii). This concludes the proof of Proposition~\ref{reacgoodsite}. \qed

\noindent\textit{Proof of Proposition \ref{UCLTprop1}.}
Let us prove (\ref{UCLTrefor}). Consider $x\in [-H\sqrt{n},H\sqrt{n}]^d$. 
We start by writing
\begin{align}
\label{UFin}
\Big| {\mathtt E}_{\theta_x\omega}[F(Z^{n})]-E[F(W)]\Big|
&\leq \Big|{\mathtt E}_{\theta_x\omega}\Big(\!F(Z^{n})
- {\mathtt E}_{\theta_{X_{\tau_{\G}}} \omega}[F(Z^{n})]\Big)\Big|
+ \Big|{\mathtt E}_{\theta_x\omega}
\Big({\mathtt E}_{\theta_{X_{\tau_{\G}}}
\omega}[F(Z^{n})]-E[F(W)]\Big)\Big|\nonumber\\
&:=U+V.
\end{align}
First, taking $\eps \leq \frac{\tilde{\eps}}{2}$ we obtain $V\leq \tilde{\eps}/2$
by definition of a $(\eps,n)$-good site.
It remains to treat the first term of the right-hand side
of~(\ref{UFin}). Denote $X':=X-x$. Now, observe that by the Markov property
\begin{align}
\label{ZION-1}
U
= \Big|{\mathtt E}_{\theta_x\omega}\Big(F(Z^{n})-{\mathtt
E}_{\theta_{X'_{\tau_{\G}}}( \theta_x
\omega)}[F(Z^{n})]\Big)\Big| \leq {\mathtt E}_{\theta_x \omega}
\Big|F\circ Z^{n}-F\circ \theta_{n^{-1}\tau_{\G}}(
Z^{n}-n^{-1/2}X'_{\tau_\G})\Big|.
\end{align}
We are going to show that for~$n$ sufficiently large we 
have uniformly in $x\in [-H\sqrt{n},H\sqrt{n}]^d$,
\[
{\mathtt E}_{\theta_x\omega}\Big|F\circ Z^{n}
 -F\circ \theta_{n^{-1}\tau_{\G}}(
Z^{n}-n^{-1/2}X'_{\tau_\G})\Big|\leq \frac{\tilde{\eps}}{2}
\]
for small enough $\eps$.
Let $M^{n}:=Z^{n}-n^{-1/2}X'_{\tau_\G}$.
Since~$F$ is uniformly continuous, we can choose $\eta>0$ in such
a way that if $\dd(f,g)\leq\eta$ then $|F(f)-F(g)|\leq
\frac{\tilde{\eps}}{4}$. Then, we have
\begin{align}
\label{ZION0}
{\mathtt E}_{\theta_x\omega} \Big|F\circ
Z^{n}-F\circ\theta_{n^{-1}\tau_{\G}}M^{n} \Big|
&={\mathtt E}_{\theta_x\omega}
\Big[ \Big|F\circ
Z^{n}-F\circ\theta_{n^{-1}\tau_{\G}}M^{n} 
\Big|\1{\dd(Z^{n},\theta_{n^{-1}\tau_{\G}}M^{n})\leq 
\eta}\Big]\nonumber\\
&\phantom{**}
 +{\mathtt E}_{\theta_x\omega}\Big[\Big|F\circ
Z^{n}-F\circ\theta_{n^{-1}\tau_{\G}}M^{n} 
\Big|\1{\dd(Z^{n},\theta_{n^{-1}\tau_{\G}}M^{n})> 
\eta}\Big]\nonumber\\
&\leq \frac{\tilde{\eps}}{4}+2\|F\|_{\infty}
{\mathtt
P}_{\theta_x\omega}\Big[\dd(Z^{n},
\theta_{n^{-1}\tau_{\G}}M^{n})>\eta\Big].
\end{align}
Since $h_{\eps}\leq1$, we have
\begin{align}
\label{ZION}
{\mathtt P}_{\theta_x\omega}
\Big[\dd(Z^{n},\theta_{n^{-1}\tau_{\G}}M^{n})>
\eta\Big]
&\leq {\mathtt P}_{\theta_x\omega}\Big[\dd(Z^{n},
\theta_{n^{-1}\tau_{\G}}M^{n})>\eta,
\tau_{\mathcal{G}}\leq hn\Big]+{\mathtt
P}_{\theta_x\omega}[\tau_{\mathcal{G}}> h_{\eps}n]\nonumber\\
&\leq {\mathtt P}_{\theta_x\omega}
\Big[\sup_{t\in [0, n^{-1}\tau_\G]}\|
Z^{n}-\theta_{n^{-1}\tau_{\G}}M^{n}\|>\frac{\eta}{2},
\tau_{\mathcal{G}}\leq h_{\eps}n\Big]\nonumber\\
&\phantom{**}+{\mathtt P}_{\theta_x\omega}
\Big[\dd(\theta_{n^{-1}\tau_{\G}}
Z^{n},\theta^2_{n^{-1}\tau_{\G}}M^{n})>\frac{\eta}{2},
\tau_{\mathcal{G}}\leq h_{\eps}n\Big]+{\mathtt
P}_{\theta_x\omega}[\tau_{\mathcal{G}}> h_{\eps}n].
\end{align}
Let $\mathcal{F}_{\tau_\G}$ be the $\sigma$-field generated by~$X$
until time $\tau_\G$. We first decompose the first term of the
right-hand side of~(\ref{ZION}) in the following way:
\begin{align}
\label{ZION1}
\lefteqn{
{\mathtt P}_{\theta_x\omega}
\Big[\sup_{t\in [0, n^{-1}\tau_\G]}\|
Z^{n}-\theta_{n^{-1}\tau_{\G}}M^{n}\|>\frac{\eta}{2},
\tau_{\mathcal{G}}\leq
h_{\eps}n\Big]
}\phantom{***********}\nonumber\\
&\leq {\mathtt P}_{\theta_x\omega}
\Big[\sup_{t\in [0, n^{-1}\tau_\G]}\|
Z^{n}\|>\frac{\eta}{4}\Big]+{\mathtt
P}_{\theta_x\omega}\Big[\sup_{t\in [0,
h_{\eps}]}\|\theta_{n^{-1}\tau_{\G}}M^{n}\|>
\frac{\eta}{4}\Big]\nonumber\\
&= {\mathtt P}_{\theta_x\omega}
\Big[\sup_{t\in [0, n^{-1}\tau_\G]}\|
Z^{n}\|>\frac{\eta}{4}\Big]+{\mathtt
E}_{\theta_x\omega}\Big(\mathtt{P}_{\theta_x\omega}\Big[\sup_{t\in
[0, h_{\eps}]}\|\theta_{n^{-1}\tau_{\G}} M^{n}\|>
\frac{\eta}{4}\mid \mathcal{F}_{\tau_\G}
\Big]\Big)\nonumber\\
&= {\mathtt P}_{\theta_x\omega}
\Big[\sup_{t\in [0, n^{-1}\tau_\G]}\|
Z^{n}\|>\frac{\eta}{4}\Big]+{\mathtt
E}_{\theta_x\omega}\Big(\mathtt{P}_{\theta_{X_{\tau_\G}}\omega}\Big[
\sup_{t\in [0, h_{\eps}]}\|Z^{n}\|>\frac{\eta}{4}\Big]\Big).
\end{align}
We now deal with the second term of the right-hand side
of~(\ref{ZION}):
\begin{align}
\label{ZION2}
\lefteqn{{\mathtt P}_{\theta_x\omega}
\Big[\dd(\theta_{n^{-1}\tau_{\G}}
Z^{n},\theta^2_{n^{-1}\tau_{\G}}M^{n})>
\frac{\eta}{2},\tau_{\mathcal{G}}\leq
h_{\eps}n\Big]}\phantom{*******}\nonumber\\
&\leq {\mathtt P}_{\theta_x\omega}\Big[\|X'_{\tau_\G}\|>\frac{\eta}{4}n\Big]+
{\mathtt
P}_{\theta_x\omega}\Big[\dd(\theta_{n^{-1}\tau_{\G}}
M^{n},\theta^2_{n^{-1}\tau_{\G}}M^{n})>
\frac{\eta}{4},\tau_{\mathcal{G}}\leq h_{\eps}n\Big]\nonumber\\
&\leq  {\mathtt P}_{\theta_x\omega}\Big[\sup_{t\in [0,
n^{-1}\tau_\G]}\| Z^{n}\|>\frac{\eta}{4}\Big]
+{\mathtt E}_{\theta_x\omega}
\Big(\1{\tau_{\mathcal{G}}\leq h_{\eps}n}{\mathtt
P}_{\theta_x\omega}\Big[\dd(\theta_{n^{-1}\tau_{\G}}
M^{n},\theta^2_{n^{-1}\tau_{\G}}M^{n})>
\frac{\eta}{4}\mid \mathcal{F}_{\tau_\G}\Big]\Big)\nonumber\\
&= {\mathtt P}_{\theta_x\omega}
\Big[\sup_{t\in [0, n^{-1}\tau_\G]}\|
Z^{n}\|>\frac{\eta}{4}\Big]
+{\mathtt E}_{\theta_x\omega}
\Big(\1{\tau_{\mathcal{G}}\leq
h_{\eps}n}\mathtt{P}_{\theta_{X_{\tau_\G}}\omega}
\Big[\dd(Z^{n},\theta_{n^{-1}\tau_{\G}}Z^{n})
>\frac{\eta}{4}\Big]\Big).
\end{align}
Combining~(\ref{ZION}), (\ref{ZION1}) and~(\ref{ZION2}), we obtain
\begin{align}
\label{ZION3}
{\mathtt P}_{\theta_x\omega}\Big[\dd(Z^{n},
\theta_{n^{-1}\tau_{\G}}M^{n})>\eta\Big]
&\leq {\mathtt P}_{\theta_x\omega}[\tau_{\mathcal{G}}> h_{\eps}n]
+2{\mathtt P}_{\theta_x\omega}\Big[\sup_{t\in [0, n^{-1}\tau_\G]}\|
Z^{n}\|>\frac{\eta}{4}\Big]
 \nonumber\\
 &\phantom{**}
+{\mathtt E}_{\theta_x\omega}
\Big(\mathtt{P}_{\theta_{X_{\tau_\G}}\omega}\Big[\sup_{t\in [0,
h_{\eps}]}\|Z^{n}\|>\frac{\eta}{4}\Big]
 \nonumber\\
 &\phantom{*******}
+\1{\tau_{\mathcal{G}}\leq h_{\eps}n}
\mathtt{P}_{\theta_{X_{\tau_\G}}\omega}
\Big[\dd(Z^{n},\theta_{n^{-1}\tau_{\G}}Z^{n})
>\frac{\eta}{4}\Big]\Big).
\end{align}
On one hand, by definition of 
a $(\eps,n)$-good point, choosing
small enough $\eps>0$,
 we have uniformly in $x\in
[-H\sqrt{n},H\sqrt{n}]^d$,
\begin{align}
\label{XION}
{\mathtt E}_{\theta_x\omega}
\Big(\mathtt{P}_{\theta_{X_{\tau_\G}}\omega}\Big[\sup_{t\in [0,
h_{\eps}]}\|Z^{n}\|>\frac{\eta}{4}\Big]
&
+\1{\tau_{\mathcal{G}}\leq
h_{\eps}n}\mathtt{P}_{\theta_{X_{\tau_\G}}\omega}
\Big[\dd(Z^{n},\theta_{n^{-1}\tau_{\G}}Z^{n})
>\frac{\eta}{4}\Big]\Big)\leq
\frac{\tilde{\eps}}{32\|F\|_{\infty}}
\end{align}
for all sufficiently large~$n$.
On the other hand, by Proposition~\ref{reacgoodsite}, for sufficiently small~$\eps$,  we have uniformly in $x\in [-H\sqrt{n},H\sqrt{n}]^d$,
\begin{equation}
\label{XION1}
{\mathtt P}_{\theta_x\omega}[\tau_{\mathcal{G}}> h_{\eps}n]
\leq \frac{\tilde{\eps}}{32\|F\|_{\infty}}
\phantom{**}\mbox{and}\phantom{**}
{\mathtt P}_{\theta_x\omega}\Big[\sup_{t\in [0, n^{-1}\tau_\G]}\|
 Z^{n}\|>\frac{\eta}{4}\Big]\leq
\frac{\tilde{\eps}}{32\|F\|_{\infty}}
\end{equation}
for sufficiently large $n$.
Combining~(\ref{XION}), (\ref{XION1})
with~(\ref{ZION3}), (\ref{ZION}), (\ref{ZION0}) and~(\ref{ZION-1}),
we have
$U\leq \tilde{\eps}/2$.
Together with $V\leq \tilde{\eps}/2$, 
this leads to the desired result.
\qed

Denote by $\C_b(C(\R_+),\R)$ 
the space of bounded continuous functionals from $C(\R_+)$ into~$\R$ and by $\mathcal{B}$ the Borel $\sig$-field on $C(\R_+)$. 
The next step is the following proposition, 
its proof follows essentially the proof of
Theorem 2.1 of~\cite{Billing} (cf. also Proposition 3.7 of~\cite{GP}).
\begin{prop}
\label{portmant}
The first statement implies the second one:

\begin{itemize}
\item[(i)] for any $F\in \C^u_b(C(\R_+),\R)$, we have $\IP$-a.s., 
\[
\lim_{n \to \infty} \sup_{x\in [-H\sqrt{n},H\sqrt{n}]^d}\Big| 
 {\mathtt E}_{\theta_x \omega}[F(Z^{n})]-E[F(W)]\Big|=0;
 \]
\item[(ii)] for any open set $G$, we have $\IP$-a.s., 
\[
\liminf_{n \to \infty} \inf_{x\in [-H\sqrt{n},H\sqrt{n}]^d}{\mathtt
P}_{\theta_x \omega}[Z^{n}\in G]\geq P[W\in G].
 \]
\end{itemize}
\end{prop}

 Finally, we have Proposition \ref{portlast}, which is similar to 
Proposition 3.8 of~\cite{GP}.
\begin{prop}
\label{portlast}
The following statements are equivalent:
\begin{itemize}
\item[(i)] we have $\IP$-a.s., for every open set~$G$,
\[
\liminf_{n \to \infty} 
\inf_{x\in [-H\sqrt{n},H\sqrt{n}]^d}{\mathtt P}_{\theta_x
\omega}[Z^{n}\in G]\geq P[W\in G];
 \]
\item[(ii)] for every open set~$G$, we have $\IP$-a.s., 
\[
\liminf_{n \to \infty} 
\inf_{x\in [-H\sqrt{n},H\sqrt{n}]^d}{\mathtt P}_{\theta_x
\omega}[Z^{n}\in G]\geq P[W\in G].
 \]
\end{itemize}
\end{prop}

\noindent
\textit{Proof.} (i)~$\Rightarrow$~(ii) is trivial. Let us
show that (ii)~$\Rightarrow$~(i). Suppose that there exists a
countable family $\mathcal{H}$ of open sets such that for every open
set~$G$ there exists a sequence $(O_n)_{n=1,2,\ldots} \subset
\mathcal{H}$ such that ${\bf 1}_{O_n} \uparrow {\bf 1}_{G}$ pointwise
as $n\to \infty$. By~(ii), since the family $\mathcal{H}$ is
countable we would have, $\IP$-a.s., for all $O\in \mathcal{H}$, 
\begin{equation}
\liminf_{n \to \infty} \inf_{x\in [-H\sqrt{n},H\sqrt{n}]}
{\mathtt P}_{\theta_x \omega}[Z^{n}\in O]\geq P[W\in O].
\end{equation}
Then, the same kind of reasoning as that used in the proof of
Proposition~\ref{portmant} to prove (i)~$\Rightarrow$~(ii) would
provide the desired result. The fact that $\mathcal{H}$ exists,
 follows from the fact that the space $C(\R_+)$ is second-countable.
\qed

\noindent 
\textit{Proof of Theorem~\ref{Theouni}.}
One can check that it is straightforward (using the same 
arguments as in the proof of Proposition~3.3) to deduce that (i), (ii), (iii) and (v) of Theorem \ref{Theouni} are equivalent to statement~(i) of Proposition~\ref{portlast}. 
That is, one can prove the equivalence of items (i)-(v) of Theorem \ref{Theouni}. To conclude the proof of Theorem~\ref{Theouni}, it remains to show that (ii) of Proposition~\ref{portlast} holds.
By Proposition~\ref{portmant}, (ii) of Proposition~\ref{portlast} is equivalent to (i) of Proposition~\ref{portmant}. Since by Proposition~\ref{UCLTprop1}, (i) of Proposition~\ref{portmant} holds, the proof of Theorem \ref{Theouni} is complete.
\qed


\section{Proof of Theorem~\ref{Theocond}}
\label{s_proof_Theocond}
For the sake of brevity, let us denote in this section, the process $DZ^n$ (resp.\ $DX$) by $\Zc$ (resp. $\Xc$). We also recall that $W^{(d)}=(W_1,\dots, W_d)$ is a $d$-dimensional standard Brownian motion.
In order to prove Theorem~\ref{Theocond}, we first show
convergence of the finite-dimensional distributions and then, in
Section~\ref{Tighteunesse}, 
we prove the tightness of the sequence
$(\Po[\Zc^n\in \cdot \mid \Lambda_n])_{n\geq 1}$. 
For $\eps\in (0,1)$, we recall that  $N:=\lfloor \eps \sqrt{n}
\rfloor$.
In this section for any set $F\subset \R^d$ we denote 
\[
\beta_F=\inf\{n\geq 0: \Xc(n)\in F\}\phantom{**} \mbox{and}\phantom{**} \beta^+_F=\inf\{n\geq 1: \Xc(n)\in F\}.
\]

\subsection{Convergence of finite-dimensional distributions}

First, let us prove
\begin{prop}
\label{propcltord}
We have $\IP$-a.s.,
\begin{equation}
\label{propordclt}
\lim_{n\to \infty}\Po[\Zc^n_1(1) > u_1, \dots, \Zc^n_d(1)>u_d \mid \Lambda_n]
=\exp(-u_1^2/2) \prod_{i=2}^{d}\int_{u_i}^{\infty}
\frac{e^{-\frac{t^2}{2}}}{\sqrt{2\pi}}dt,
\end{equation}
for all $u=(u_1,\dots, u_d)\in \R_+\times \R^{d-1}$.
\end{prop}

\noindent
\textit{Proof.}
First, we introduce some notations. Let
\[
\mathcal{D}_{u}=\{x\in \R^d: x_1>u_1,\dots, x_d>u_d\}
\]
and 
\[
R_{\eps, n}=\{x\in \R^d: x_1=N, x_i\in [-\lf \eps^{-1/2}N \rf,
 \lf \eps^{-1/2} N\rf], i\in[\![2,d]\!]\}.
\]
Let us denote $\mathcal{R}_{\eps, n}=DR_{\eps, n}$, we also define the event 
$A_{0\to R}=\{\beta_{\Rc_{\eps, n}}<\beta_{\{0\}_1}^+\}$. We start by
bounding the term $\Po[\Zc^n(1) \in \mathcal{D}_{u}\mid \Lambda_n]$ from above.
Fix $\eps\in (0, u_1\wedge 1)$ and consider the 
following decomposition
\begin{align}
\label{FI1}
\Po[\Zc^n(1) \in \mathcal{D}_{u}\mid
\Lambda_n]
&\leq \frac{1}{\Po[\Lambda_n]}\Big(\Po[\Zc^n(1) \in \mathcal{D}_{u},A_{0\to R},\Lambda_n]
+\Po[A^c_{0\to R},\Lambda_n]\Big)\nonumber\\
&=\frac{1}{\Po[\Lambda_n]}\Big(\Po[\Zc^n(1) \in \mathcal{D}_{u},A_{0\to
R},\Lambda_n, \beta_{\Rc_{\eps, n}}\leq \eps^{1/2} n]\nonumber\\
&\phantom{*******}+\Po[\Zc^n(1) \in \mathcal{D}_{u},A_{0\to R},\Lambda_n, 
\beta_{\Rc_{\eps, n}}> \eps^{1/2} n]   \Big)
+\Po[A^c_{0\to R}\mid \Lambda_n]\nonumber\\
&\leq (\Po[\Lambda_n])^{-1}\Po[\Zc^n(1) \in \mathcal{D}_{u},A_{0\to R},\Lambda_n, 
\beta_{\Rc_{\eps, n}}\leq \eps^{1/2} n]
+\Po[\beta_{\Rc_{\eps, n}}> \eps^{1/2} n\mid \Lambda_n]\nonumber\\
&\phantom{**}+\Po[A^c_{0\to R}\mid \Lambda_n].
\end{align}
Since $\eps^{1/2}\in (0,1)$, we have 
\begin{equation}
\label{decompo1}
\Po[A^c_{0\to R}\mid \Lambda_n]
= \Po[\beta_{\Rc_{\eps, n}}>\beta^+_{\{0\}_1}\mid \Lambda_n]
\leq \Po[\beta_{\Rc_{\eps, n}}>n\mid \Lambda_n]
\leq \Po[\beta_{\Rc_{\eps, n}}>\eps^{1/2} n\mid \Lambda_n].
\end{equation}
Then, using the Markov property at time $\beta_{\Rc_{\eps,n}}$ we deduce
\begin{align}
\label{FI2}
\frac{1}{\Po[\Lambda_n]}\Po[\Zc^n(1) \in \Dc_{u},A_{0\to R},\Lambda_n,
 \beta_{\Rc_{\eps, n}}\leq \eps^{1/2} n]
&\leq \frac{\Po[A_{0\to R}]}{\Po[\Lambda_n]}\max_{y\in R_{\eps, n}}\max_{j\leq 
\lf \eps^{1/2} n\rf}\Po^y\Big[\frac{\Xc(n-j)}{\sqrt{n}}\in \Dc_{u},\Lambda_{n-j}\Big].
\end{align}
Again, using the Markov property at time $\beta_{\Rc_{\eps,n}}$ we obtain
\begin{align}
\label{FI3}
\frac{\Po[\Lambda_n]}{\Po[A_{0\to R}]}
\geq \min_{y\in R_{\eps, n}}\Po^y[\Lambda_{n}].
\end{align}
Combining (\ref{FI1}), (\ref{decompo1}), (\ref{FI2}) and (\ref{FI3}) we obtain
\begin{align}
\Po[\Zc^n(1)\in \Dc_{u}\mid \Lambda_n]\leq \frac{\max_{y\in R_{\eps, n}}
\max_{j\leq  \lf \eps^{1/2} n\rf}\Po^y[\Xc(n-j)\in \Dc_{u}\sqrt{n},\Lambda_{n-j}]}
{\min_{y\in R_{\eps, n}}\Po^y[\Lambda_{n}]}+2\Po[\beta_{\Rc_{\eps, n}}> 
 \eps^{1/2} n\mid \Lambda_n].
\end{align}
Now, to bound the term 
$\Po[\beta_{\Rc_{\eps, n}}> \eps^{1/2} n\mid \Lambda_n]$ from above we notice that
\begin{align}
\Po[\beta_{\Rc_{\eps, n}}>  \eps^{1/2} n\mid \Lambda_n]
&=\Po[\tau_{R_{\eps, n}}>  \eps^{1/2} n\mid \Lambda_n]\nonumber\\
&\leq \Po\Big[\max_{i\in [\![2,d]\!]}\sup_{j\leq \tau_{\{N\}_1}}|X_i(j)|>
\eps^{-1/2}N\mid \Lambda_n\Big]+\Po[\tau_{\{N\}_1}> \eps^{1/2} n\mid \Lambda_n].
\end{align}
By Lemmas~\ref{Lem2} and~\ref{Lem3} we have
\begin{equation}
\label{URF}
\limsup_{n \to \infty}\Po[\beta_{\Rc_{\eps, n}}> 
 \eps^{1/2} n\mid \Lambda_n]\leq f(\eps)+g(\eps).
\end{equation}
By definition of $\Zc^n$, we have $\Po^y[\Lambda_{n}]=\Po^{y}\Big[\Zc^n_1(1)>0, t\in [0,1]\Big]$.
Thus, from Theorem~\ref{Theouni} we obtain, recalling that $W_1$ is the first component of $W^{(d)}$,
\begin{align}
\label{lamb1}
\lim_{n\to \infty}\min_{y\in R_{\eps, n}}\Po^{y}\Big[\Zc^n_1(t)>0, t\in [0,1]\Big]
&=P^{\eps \sig_1}\Big[\min_{0\leq t\leq 1} W_1(t)>0\Big]
=P\big[|W_1(1)|<\eps \sig_1\big]\nonumber\\
&=\frac{2\eps \sig_1}{ \sqrt{2\pi}}+o(\eps) 
\end{align}
as $\eps\to 0$, where $P^x$ is law of $W^{(d)}$ starting at $x$ and $\sig_1:=D{\bf{e}}_1\cdot {\bf{e}}_1>0$ (cf.\ Section \ref{s_intro}).
Now, let us treat the term 
$$\max_{y\in R_{\eps, n}}\max_{j\leq \lf \eps^{1/2}
n\rf}\Po^y[\Xc(n-j)\in \Dc_{u}\sqrt{n},\Lambda_{n-j}].$$ Fix $\delta'>0$ and let $\sgn(x)=-1$ if $x\leq 0$ and $1$ if $x>0$. Denote by 
\begin{equation}
\label{def_Ti}
U_i:=\Big\{ \Xc_i(n-\lfloor\eps^{1/2} n \rfloor) >  (u_i-\sgn(u_i)\delta')
\sqrt{n}\Big\}
\end{equation}
and 
\begin{equation}
\label{def_Vi}
V_i:=\Big\{\max_{j\leq \lf \eps^{1/2} n\rf}|\Xc_i(n-\lfloor \eps^{1/2} n \rfloor)
-\Xc_i(n-j)|\geq \delta'
\sqrt{n}\Big\}
\end{equation}
for $i=1,\dots,d$. Observe that we have for $y\in R_{\eps, n}$ and
$j\leq \lf\eps^{1/2} n\rf$
\begin{align*}
\Po^y[\Xc(n-j)\in \Dc_{u}\sqrt{n},\Lambda_{n-j}]
 &\leq \Po^{y}\Big[\bigcap_{i=1}^{d}(U_i\cup V_i)\cap  
\Lambda_{n-\lf \eps^{1/2} n\rf}\Big].
\end{align*}
Let us consider the set $\mathcal{I}=\{U_1,\dots,U_d, V_1,\dots, V_d\}$ 
and denote by $\mathcal{J}$ the set formed by all intersections of~$d$
 distinct elements of $\mathcal{I}$: $\mathcal{J}$ contains
 ${\binom{2d}{d}}$ elements. Let us denote by $J_1,\dots, J_{{\binom{2d}{d}}}$ 
all the elements of $\mathcal{J}$.
Therefore, we obtain  
\begin{align}
\label{DEC1}
\max_{j \leq \lfloor \eps^{1/2} n \rfloor}
\Po^y[\Xc(n-j)\in \Dc_{u}\sqrt{n},\Lambda_{n-j}]\leq
\sum_{i\leq{\binom{2d}{d}}}\Po^y\Big[J_i, \Lambda_{n-\lf \eps^{1/2} n \rf}\Big].
\end{align}
Let us treat the term 
$\Po^y[\cap_{i=1}^{d}U_i, \Lambda_{n-\lf \eps^{1/2} n \rf}]$. 
We have by definition of $\Zc^n$
\begin{align*}
\Po^y\Big[\bigcap_{i=1}^{d}U_i, \Lambda_{n-\lf \eps^{1/2} n \rf}\Big]
\leq\Po^{y}\Big[\bigcap_{i=1}^{d}\Big\{\Zc_i^{n-\lf \eps^{1/2} n\rf}(1)
>(u_i-\sgn(u_i)\delta')\Big\},\Zc_1^{n-\lf \eps^{1/2} n\rf}(t)>0, t\in [0,1]\Big].
\end{align*}
By Theorem \ref{Theouni} we deduce 
\begin{align}
\label{Autre1}
\limsup_{n \to \infty}\max_{y\in R_{\eps,n}}
\Po^y\Big[\bigcap_{i=1}^{d}U_i, \Lambda_{n-\lf \eps^{1/2} n
\rf}\Big]
&\leq P^{\frac{\eps\sig_1}{\sqrt{1-\eps^{1/2}}}}\Big[W_1(1)> 
(u_1-\sgn(u_1)\delta'), \min_{0\leq t\leq
1}W_1(t)>0\Big]\nonumber\\
&\phantom{**}\times \prod_{i=2}^{d}
P^{\frac{\gam_1\eps^{1/2}}{\sqrt{1-\eps^{1/2}}}}
[W_i(1)>(u_i-\sgn(u_i)\delta')]
\end{align}
for some constant $\gam_1$.
 Abbreviate $\eps':=\sig_1\eps (1-\eps^{1/2})^{-1/2}$ and let us compute
the first term of the right-hand side of~(\ref{Autre1}) 
for sufficiently small~$\eps$.
By the reflection principle for the Brownian motion, we have
\begin{align*}
\lefteqn{
P^{\eps'}\Big[W_1(1) > (u_1-\sgn(u_1)\delta'),\min_{0\leq t\leq 1}W_1(t)>0\Big]
}\phantom{***********}\nonumber\\
 &=P^{\eps'}\Big[W_1(1) > (u_1-\sgn(u_1)\delta')\Big]-P^{\eps'}\Big[W_1(1) <- (u_1-\sgn(u_1)\delta')\Big]\nonumber\\
 &=P\Big[W_1(1)> (u_1-\sgn(u_1)\delta')-\eps'\Big]-P\Big[W_1(1)
<-(u_1-\sgn(u_1)\delta')-\eps')\Big]\\
&=\frac{1}{\sqrt{2\pi}}\int_{(u_1-\sgn(u_1)\delta')-\eps'}^
{(u_1-\sgn(u_1)\delta')+\eps'}e^{-\frac{x^2}{2}}dx.
\end{align*}
Therefore, we obtain, as $\eps\to 0$
\begin{align}
\label{SpL2}
\lefteqn{\limsup_{n \to \infty}\max_{y\in R_{\eps,n}}
\Po^y\Big[\bigcap_{i=1}^{d}U_i, \Lambda_{n-\lf \eps^{1/2} n
\rf}\Big]}\phantom{****************}\nonumber\\
&\leq\Big(\frac{2\eps\sig_1 e^{-\frac{(u_1-\sgn(u_1)\delta')^2}{2}}}{
\sqrt{2\pi(1-\eps^{1/2})}}+o(\eps)\Big)\prod_{i=2}^{d}\int_{
(u_i-\sgn(u_i)\delta')-\frac{\gam_1\eps^{1/2}}{\sqrt{1-\eps^{1/2}
}}}^{\infty}\frac{e^{-\frac{t^2}{2}}}{\sqrt{2\pi}}dt.
\end{align}

The other terms $\Po^y[J_i, \Lambda_{n-\lf \eps^{1/2} n \rf}]$ necessarily 
contain a term $V_j$ for some $j\in [\![1,d]\!]$. Thus, we have for 
$J_i\neq \cap_{i=1}^{d}U_i$,
\begin{equation}
\label{TPX}
\Po^y[J_i, \Lambda_{n-\lf \eps^{1/2} n \rf}]\leq \sum_{j=1}^{d}\Po^y[V_j].
\end{equation}
Let us bound the terms $\limsup_{n\to \infty} \max_{y\in R_{\eps,n}}\Po^y[V_j]$
 for $j\in [\![1,d]\!]$. We start by writing
\begin{align*}
\Po^y[V_j]
&=\Po^{y}\Big[\max_{i \leq \lfloor \eps^{1/2} n
\rfloor}|\Xc_j(n-\lfloor \eps^{1/2} n \rfloor)-\Xc_j(n-i)|\geq \delta' \sqrt{n}\Big]\nonumber\\
 &=\Po^{y}\Big[\max_{n-\lfloor \eps^{1/2} n \rfloor\leq k\leq
n}\Big|\Xc_j(k)-\Xc_j(n-\lfloor \eps^{1/2} n \rfloor)\Big|\geq \delta' 
\sqrt{n}\Big]\nonumber\\
 &\leq \Po^{y}\Big[\max_{1-\eps^{1/2} \leq t\leq
1}\Big(\Zc_j^n(t)-\min_{1-\eps^{1/2}\leq s\leq t}\Zc_j^n(s)\Big) \geq \delta' \Big] \nonumber\\
&\phantom{**}
+\Po^{y}\Big[\min_{1-\eps^{1/2} \leq
t\leq 1}\Big(\Zc_j^n(t)-\max_{1-\eps^{1/2}\leq s\leq t}\Zc_j^n(s)\Big) 
\leq-\delta' \Big].
\end{align*}
By Theorem~\ref{Theouni}, we obtain
\begin{align}
\label{Convfraca1}
 \lefteqn{ 
\lim_{n\to \infty}\max_{y\in R_{\eps,n}}\Po^{y}\Big[\max_{1-\eps^{1/2} 
 \leq t\leq 1}\Big(\Zc_j^n(t)-\min_{1-\eps^{1/2} \leq
s\leq t}\Zc_j^n(s)\Big) \geq \delta' \Big]}\phantom{*******************}\nonumber\\
&=P\Big[\max_{1-\eps^{1/2}  \leq t\leq
1}\Big(W_j(t)-\min_{1-\eps^{1/2} \leq s\leq t}W_j(s)\Big) \geq \delta'
\Big]
\end{align}
and 
\begin{align}
\label{Convfraca2}
\lefteqn{
\lim_{n\to \infty}\max_{y\in R_{\eps,n}}\Po^{y}\Big[\min_{1-\eps^{1/2} 
 \leq t\leq 1}\Big(\Zc_j^n(t)-\max_{1-\eps^{1/2} \leq
s\leq t}\Zc_j^n(s)\Big) \leq-\delta' \Big]}\phantom{*******************}\nonumber\\
 &=P\Big[\min_{1-\eps^{1/2}  \leq t\leq
1}\Big(W_j(t)-\max_{1-\eps^{1/2} \leq s\leq t}W_j(s)\Big) \leq-\delta' \Big].
\end{align}
 Observe that the right-hand sides of~(\ref{Convfraca1}) 
and~(\ref{Convfraca2}) are equal since $(-W_j)$ is a Brownian motion.
Thus, let us compute for example the right-hand side term of~(\ref{Convfraca1}).
 By L\'evy's Theorem (cf.~\cite{RY}, Chapter~VI, Theorem~2.3), we
have
\begin{equation*}
P\Big[\max_{0  \leq t\leq \eps^{1/2}}\Big(W_j(t)-\min_{0 \leq s\leq t}W_j(s)\Big)
 \geq \delta'\Big]
=P\Big[\max_{0\leq t\leq
\eps^{1/2}}|W_j(t)|\geq \delta' \Big].
\end{equation*}
Then, 
\begin{align*}
P\Big[\max_{0\leq t\leq
\eps^{1/2}}|W_j(t)|\geq \delta' \Big]
&
 \leq 2P\Big[\max_{0\leq t\leq\eps^{1/2}}W_j(t)\geq \delta'\Big]= 4P[W_j(\eps^{1/2})\geq \delta'].
\end{align*}
 Using an estimate on the tail of the Gaussian law
(cf.~\cite{PerMot}, Appendix~B, Lemma~12.9) we obtain
\begin{equation*}
 P\Big[\max_{0\leq t\leq \eps^{1/2}}|W_j(t)|\geq \delta' \Big]\leq
\frac{4\eps^{1/4}}{\delta'
\sqrt{2\pi}}\exp\Big\{-\frac{(\delta')^2}{2\eps^{1/2}}\Big\}.
\end{equation*}
We finally obtain
\begin{equation}
\label{Autre2}
\limsup_{n \to \infty}\max_{y \in R_{\eps, n}}\sum_{i=1}^{d}\Po^y[V_i]
\leq \frac{8d\eps^{1/4}}{\delta'\sqrt{2\pi}}\exp\Big\{- \frac{(\delta')^2}
{2\eps^{1/2} } \Big\}.
\end{equation}
To sum up, combining 
(\ref{lamb1}), (\ref{DEC1}), (\ref{SpL2}), (\ref{TPX}), and~(\ref{Autre2}), we
have $\IP$-a.s.
\begin{align}
\label{JUGD1}
\lefteqn{\limsup_{n \to \infty} \Po[\Zc^n(1)\in \Dc_{u}\mid
\Lambda_n]}\phantom{*******}\nonumber\\
&\leq \Big(\frac{2\eps\sig_1}{
\sqrt{2\pi}}+o(\eps)\Big)^{-1}\Big(\frac{2\eps\sig_1 e^{-\frac{(u_1-\sgn(u_1)\delta')^2}{2}}}{ 
\sqrt{2\pi(1-\eps^{1/2})}}+o(\eps)\Big)\prod_{i=2}^{d}\int_{
(u_i-\sgn(u_i)\delta')-\frac{\gam_1\eps^{1/2}}{\sqrt{1-\eps^{1/2}
}}}^{\infty}\frac{e^{-\frac{t^2}{2}}}{\sqrt{2\pi}}dt\nonumber\\
&\phantom{**}+{\binom{2d}{d}
}\frac{8d\eps^{1/4}}{\delta'\sqrt{2\pi}}\exp\Big\{-
\frac{(\delta')^2}{2\eps^{1/2}} \Big\}+2(f(\eps)+g(\eps)).
\end{align}


Let us now bound the term $\Po[\Zc^n(1) \in \Dc_{u}\mid
\Lambda_n]$ from below. We have by the Markov property
\begin{align}
\label{INFDEC1}
\Po[\Zc^n(1) \in \Dc_{u}\mid \Lambda_n]
\geq \frac{\Po[A_{0\to R},\beta_{\Rc_{\eps, n}}\leq \eps^{1/2} n] }
{\Po[\Lambda_n]}\min_{y\in R_{\eps, n}}\min_{j\leq \lf \eps^{1/2} n\rf}
\Po^y[\Xc(n-j)\in \Dc_{u}\sqrt{n},\Lambda_{n-j}].
\end{align}
We first decompose the term $(\Po[\Lambda_n])^{-1}\Po[A_{0\to R},
\beta_{\Rc_{\eps, n}}\leq \eps^{1/2} n]$ in the following way
\begin{align}
\label{POIT1}
 \frac{\Po[A_{0\to R},\beta_{\Rc_{\eps, n}}
\leq \eps^{1/2} n]}{\Po[\Lambda_n]}
&=\frac{\Po[A_{0\to R}]}{\Po[\Lambda_n]}-\frac{\Po[A_{0\to R},\beta_{\Rc_{\eps, n}}
>\eps^{1/2} n]}{\Po[\Lambda_n]}\nonumber\\
&=\frac{\Po[A_{0\to R}]}{\Po[\Lambda_n]}(1-\Po[\beta_{\Rc_{\eps, n}}>\eps^{1/2} n 
\mid A_{0\to R}]).
\end{align}
Then, we write 
\begin{align}
\label{TRWER4}
\Po[\beta_{\Rc_{\eps, n}}>\eps^{1/2} n \mid A_{0\to R}]
&=\frac{\Po[\beta_{\Rc_{\eps, n}}>\eps^{1/2} n, A_{0\to R}]}{\Po[A_{0\to R}]}
\leq \frac{\Po[\beta_{\Rc_{\eps, n}}>\eps^{1/2} n, \Lambda_{\eps^{1/2} n}]}
{\Po[A_{0\to R},\Lambda_{\eps^{1/2} n}]}\nonumber\\
&=\frac{\Po[\beta_{\Rc_{\eps, n}}>\eps^{1/2} n\mid \Lambda_{\eps^{1/2} n}]}
{1-\Po[A^c_{0\to R}\mid \Lambda_{\eps^{1/2} n}]}.
\end{align}
For the term $\Po[\beta_{\Rc_{\eps, n}}>\eps^{1/2} n\mid \Lambda_{\eps^{1/2} n}]$,
 we have, recalling that $N =\lf \eps \sqrt{n} \rf$,
\begin{align}
\label{INGUEE}
\Po[\beta_{\Rc_{\eps, n}}>\eps^{1/2} n\mid \Lambda_{\eps^{1/2} n}]
&=\Po[\tau_{R_{\eps, n}}>\eps^{1/2} n\mid \Lambda_{\eps^{1/2} n}]\nonumber\\
&\leq \Po\Big[\max_{i\in [\![2,d]\!]}\sup_{j\leq \tau_{\{N\}_1}}|X_i(j)|> 
\eps^{-1/2}N\mid \Lambda_{\eps^{1/2} n}\Big]+\Po[\tau_{\{N\}_1}
> \eps^{1/2} n\mid \Lambda_{\eps^{1/2} n}].
\end{align}
By Lemmas~\ref{Lem2} and~\ref{Lem3} we deduce
\begin{equation}
\label{POIT2}
\limsup_{n \to \infty}\Po[\beta_{\Rc_{\eps, n}}>\eps^{1/2} n\mid 
\Lambda_{\eps^{1/2} n}]\leq g(\eps^{3/4})+f(\eps^{3/4}).
\end{equation}
%
For the term $\Po[A^c_{0\to R}\mid \Lambda_{\eps^{1/2} n}]$, we write
\begin{align*}
\Po[A^c_{0\to R}\mid \Lambda_{\eps^{1/2} n}]
= \Po[\beta_{\Rc_{\eps, n}}>\beta_{\{0\}_1}^+\mid \Lambda_{\eps^{1/2}n}]\leq \Po[\beta_{\Rc_{\eps, n}}>\eps^{1/2} n\mid \Lambda_{\eps^{1/2} n}].
\end{align*} 
Hence, by~(\ref{POIT2}) we obtain
\begin{equation}
\label{POIT3}
\limsup_{n \to \infty}\Po[A^c_{0\to R}\mid
 \Lambda_{\eps^{1/2} n}] \leq f(\eps^{3/4})+g(\eps^{3/4}).
\end{equation} 

Going back to the term $(\Po[\Lambda_n])^{-1}\Po[A_{0\to R}]$ 
in~(\ref{POIT1}), we write
\begin{align}
\label{TRWER1}
\frac{\Po[A_{0\to R}]}{\Po[\Lambda_n]}
&=\frac{\Po[A_{0\to R}]}{\Po[\Lambda_n, A_{0\to R}]
+\Po[\Lambda_n, A^c_{0\to R}]}\nonumber\\
&=\Big(\Po[\Lambda_n \mid A_{0\to R}]+\Po[\Lambda_n, A^c_{0\to R}]
(\Po[A_{0\to R}])^{-1}\Big)^{-1}\nonumber\\
&\geq\Big(\Po[\Lambda_n \mid A_{0\to R}]+\Po[\Lambda_n, A^c_{0\to R}]
(\Po[\Lambda_n, A_{0\to R}])^{-1}\Big)^{-1}\nonumber\\
&=\Big(\Po[\Lambda_n \mid A_{0\to R}]+\Po[A^c_{0\to R}
\mid \Lambda_n](1- \Po[A^c_{0\to R}\mid \Lambda_n])^{-1}\Big)^{-1}.
\end{align}
By (\ref{URF}), we have
\begin{equation}
\label{TRWER3}
\limsup_{n\to \infty}\Po[A^c_{0\to R}\mid \Lambda_n]
\leq \limsup_{n \to \infty}\Po[\beta_{\Rc_{\eps, n}}>\eps^{1/2} n
\mid \Lambda_{ n}]
\leq f(\eps)+g(\eps).
\end{equation}
Then, we have by the Markov property
\begin{align}
\label{TRWER2}
\Po[\Lambda_n \mid A_{0\to R}]
\leq \max_{y\in R_{\eps, n}}\Po^y[\Lambda_{n-\lf \eps^{1/2} n\rf}]
+\Po[\beta_{\Rc_{\eps, n}}>\eps^{1/2} n  \mid A_{0\to R}].
\end{align}
Thus, by (\ref{TRWER1}), (\ref{TRWER2}), (\ref{TRWER4}), 
(\ref{POIT2}), (\ref{POIT3}), and~(\ref{TRWER3}), we deduce
\begin{align}
\label{POIT4}
\liminf_{n \to \infty}\frac{\Po[A_{0\to
R}]}{\Po[\Lambda_n]}
&\geq \Big( \limsup_{n \to \infty}\max_{y\in R_{\eps, n}}\Po^y[\Lambda_{n-\lf 
\eps^{1/2} n\rf}]+\frac{f(\eps^{3/4})+g(\eps^{3/4})}{1-f(\eps^{3/4})
-g(\eps^{3/4})}+\frac{f(\eps)+g(\eps)}{1-f(\eps)-g(\eps)}\Big)^{-1}.
\end{align}
Combining (\ref{INFDEC1}), (\ref{POIT1}), (\ref{POIT2}), (\ref{POIT3}), 
and~(\ref{POIT4}), we obtain $\IP$-a.s.
\begin{align}
\label{TIOFFF1}
\lefteqn{\liminf_{n \to \infty}\Po[\Zc^n(1) \in \Dc_{u}\mid \Lambda_n]}
\phantom{*******}\nonumber\\
&\geq \Big(\limsup_{n \to \infty} \max_{y\in R_{\eps, n}}\Po^y[\Lambda_{n-\lf 
\eps^{1/2} n\rf}]+\frac{f(\eps^{3/4})+g(\eps^{3/4})}{1-f(\eps^{3/4})
-g(\eps^{3/4})}+\frac{f(\eps)+g(\eps)}{1-f(\eps)-g(\eps)}\Big)^{-1}\nonumber\\
&\phantom{**}\times \Big(1-\frac{f(\eps^{3/4})+g(\eps^{3/4})}{1-f(\eps^{3/4})
-g(\eps^{3/4})}\Big)\liminf_{n \to \infty}\min_{y\in R_{\eps, n}}
\min_{j\leq \eps^{1/2} n}\Po^y[\Xc(n-j)\in \Dc_{u}\sqrt{n},\Lambda_{n-j}].
\end{align}
Analogously to (\ref{lamb1}) we have
\begin{align}
\label{TIOFFF2}
\lim_{n\to \infty}\max_{y\in R_{\eps, n}}\Po^y[\Lambda_{n-\lf \eps^{1/2} n\rf}]
&=\frac{2\eps \sig_1}{ \sqrt{2\pi(1-\eps^{1/2})}}+o(\eps).
\end{align}
At this point, let us introduce more notations. Let $\delta'>0$
be the constant used in the definitions of~$V_i$ and~$U_i$ 
(cf.~\eqref{def_Ti} and~\eqref{def_Vi}) and introduce
\[
E_i=\Big\{ \Xc_i(n) > (u_i+\sgn(u_i)\delta')\sqrt{n}\Big\}\phantom{**} \mbox{and}\phantom{**} F_i=\Big\{\max_{j\leq \lf \eps^{1/2} n\rf}|\Xc_i(n)-\Xc_i(n-j)|\leq   \delta' 
\sqrt{n}\Big\}
\]
for $i\in [\![1,d]\!]$. 
Observe that for all $y\in R_{\eps,n }$ and $j\leq \lf\eps^{1/2} n\rf$ we have
\begin{align}
\Po^y[\Xc(n-j)\in \Dc_{u}\sqrt{n},\Lambda_{n-j}]
&\geq \Po^y\Big[\bigcap_{i=1}^{d}(E_i\cap F_i),\Lambda_{n}\Big] \geq \Po^y\Big[\bigcap_{i=1}^{d}E_i, \Lambda_n\Big]-\sum_{i=1}^{d}\Po^y[F_i^c].
\end{align}
By Theorem~\ref{Theouni} and similar computations as those to derive 
equations~(\ref{SpL2}) and~(\ref{Autre2}), we obtain for some constant $\gam_2$,
\begin{equation}
\label{TIOFFF3}
\lim_{n \to \infty}\min_{y\in R_{\eps,n}}
\Po^y\Big[\bigcap_{i=1}^{d}E_i, \Lambda_n\Big]
=\Big(\frac{2\eps\sig_1}{\sqrt{2\pi}}e^{-\frac{(u_1+\sgn(u_1)\delta')^2}{2}}+o(\eps)\Big)
\prod_{i=2}^{d}\int_{(u_i+\sgn(u_i)\delta')-\gam_2\eps^{1/2}}^{\infty}
\frac{e^{-\frac{t^2}{2}}}{\sqrt{2\pi}}dt
\end{equation}
 as $\eps \to 0$ and 
\begin{equation}
\label{TIOFFF4}
\limsup_{n \to \infty}\max_{y\in R_{\eps, n}}\sum_{i=1}^{d}\Po^y[F_i^c]\leq 
 \frac{8d\eps^{1/4}}{\delta' \sqrt{2\pi}}
\exp\Big\{ -\frac{(\delta')^2}{2\eps^{1/2}}\Big\}.
\end{equation}
Combining (\ref{TIOFFF1}), (\ref{TIOFFF2}), (\ref{TIOFFF3}), 
and~(\ref{TIOFFF4}), we obtain $\IP$-a.s.
\begin{align}
\label{JUGD2}
\lefteqn{\liminf_{n \to \infty}\Po[\Zc^n(1) \in \Dc_{u}\mid \Lambda_n]}
\phantom{**}\nonumber\\
&\geq \Big(\frac{2\eps\sig_1}{\sqrt{2\pi(1-\eps^{1/2})}}+o(\eps)
+\frac{f(\eps^{3/4})+g(\eps^{3/4})}{1-f(\eps^{3/4})-g(\eps^{3/4})}
+\frac{f(\eps)+g(\eps)}{1-f(\eps)-g(\eps)}\Big)^{-1}\nonumber\\
&\phantom{**}\times \Big(1-\frac{f(\eps^{3/4})+g(\eps^{3/4})}
{1-f(\eps^{3/4})-g(\eps^{3/4})}\Big)\nonumber\\
&\phantom{**}\times\Big(\Big(\frac{2\eps\sig_1}{\sqrt{2\pi}}
e^{-\frac{(u_1+\sgn(u_1)\delta')^2}{2}}+o(\eps)\Big)\prod_{i=2}^{d}\int_{(u_i+\sgn(u_i)\delta')
-\gam_2\eps^{1/2}}^{\infty}\frac{e^{-\frac{t^2}{2}}}{\sqrt{2\pi}}dt
- \frac{8d \eps^{1/4}}{\delta' \sqrt{2\pi}}
\exp\Big\{ -\frac{(\delta')^2}{2\eps^{1/2}}\Big\}\Big).
\end{align}
Finally, take $\delta'= \eps^{1/8}$ and let $\eps\to 0$ 
in~(\ref{JUGD1}) and~(\ref{JUGD2}) to prove~(\ref{propordclt}).
\qed

 The next steps in showing that the f.d.d.'s converge are 
standard and we follow~\cite{Igle} and~\cite{GP}. We start by
recalling the transition density function of the Brownian meander
(see~\cite{Igle}) from $(0,0)$ to $(t,x_1)$
\begin{equation}
\label{GRANFIN}
q(0,0;t,x_1)=t^{-3/2}x_1
\exp\Big(-\frac{{x_1}^2}{2t}\Big){\tilde N}(x_1(1-t)^{-1/2})
\end{equation}
for $x_1>0$, $0<t \leq 1$ and from $(t_1,x_1)$ to $(t_2,x_2)$
\[
 q(t_1,x_1;t_2,x_2)=g(t_2-t_1,x_1,x_2)
\frac{{\tilde N}(x_2(1-t_2)^{-1/2})}{{\tilde N}(x_1(1-t_1)^{-1/2})}
\]
for $x_1$, $x_2>0$, $0<t_1<t_2 \leq 1$, where
\[
 {\tilde N}(v)=\Big(\frac{2}{\pi}\Big)^{1/2}\int_0^v
e^{-\frac{u^2}{2}}du
\]
for $v\geq 0$ and
\[
 g(t,x_1,x_2)=(2\pi)^{-1/2}\Big(\exp\Big(-\frac{(x_2-x_1)^2}{2t}
\Big)-\exp\Big(-\frac{(x_1+x_2)^2}{2t}\Big)\Big)
\]
for $x_1$, $x_2>0$ and $0<t \leq 1$. 
\medskip

Next, we will prove the following
\begin{prop}
\label{prop_Z<x}
We have $\IP$-a.s., for $u_1> 0$, $-\infty<a_i<b_i<\infty$, $i\in[\![2,d]\!]$ 
and $0<t<1$,
\begin{equation}
\label{FFD1}
 \lim_{n\to \infty}\Po\Big[\Zc^n_1(t)\leq u_1,\bigcap_{i=2}^{d}\Big\{\Zc^n_i(t)
\in (a_i, b_i]\Big\}\mid \Lambda_n\Big]=\int_0^{u_1}
q(0,0;t,v)dv\prod_{i=2}^{d}\int_{a_i}^{b_i}\frac{e^{-\frac{v^2}{2t}}}
{\sqrt{2\pi t}}dv.
\end{equation}
\end{prop}
\textit{Proof.}
For $\eps>0$ we have
\begin{align}
\label{Finiteconvdist}
\lefteqn{\Po\Big[\Zc^n_1(n^{-1}\lfloor nt\rfloor)\leq u_1-\eps,
\bigcap_{i=2}^{d}\Big\{\Zc^n_i(n^{-1}\lfloor nt\rfloor)
\in (a_i-\eps, b_i+\eps]\Big\} \mid \Lambda_n\Big]}\phantom{*******}\nonumber\\
&\leq \Po\Big[\Zc^n_1(t)\leq u_1,\bigcap_{i=2}^{d}
\Big\{\Zc^n_i(t)\in (a_i, b_i]\Big\} \mid
\Lambda_n\Big]\nonumber\\
&\phantom{***}\leq  \Po\Big[\Zc^n_1(n^{-1}\lfloor nt\rfloor)
\leq u_1+\eps,\bigcap_{i=2}^{d}\Big\{\Zc^n_i(n^{-1}\lfloor nt\rfloor)
\in (a_i+\eps, b_i-\eps]\Big\} \mid \Lambda_n\Big].
\end{align}
for all sufficiently large~$n$. 
Now, suppose that we have for all $u_1\geq 0$, $a_i<b_i$ and $0<t<1$,
\begin{equation}
\label{Finite replace}
\lim_{n \to \infty} \Po\Big[\Zc^n_1(n^{-1}\lfloor nt\rfloor)\leq u_1,
\bigcap_{i=2}^{d}\Big\{\Zc^n_i(n^{-1}\lfloor nt\rfloor)\in (a_i, b_i]\Big\}
\mid \Lambda_n\Big]=\int_0^{u_1} q(0,0;t,v)dv
\prod_{i=2}^{d}\int_{a_i}^{b_i}\frac{e^{-\frac{v^2}{2t}}}{\sqrt{2\pi t}}dv.
\end{equation}
Combining~(\ref{Finiteconvdist}) and~(\ref{Finite replace}), we
obtain~(\ref{FFD1}) since the limit distribution $q(0,0;t,x_1)$ is
absolutely continuous. Let us denote by $l=l(t,n)$ the quantity 
$(n\lf nt\rf^{-1})^{1/2}$. 
We recall that $x$ is the vector of coordinates $(x_1,\dots,x_d)$.
Then, observe that
\begin{align}
\label{findimdist4}
\lefteqn{
\Po\Big[\Zc^n_1(n^{-1}\lfloor nt\rfloor)\leq u_1,\bigcap_{i=2}^{d}
\Big\{\Zc^n_i(n^{-1}\lfloor nt\rfloor)\in (a_i, b_i]\Big\}\mid \Lambda_n\Big]
}\phantom{**}\nonumber\\
&=\frac{1}{\Po[\Lambda_n]}\Po\Big[\Zc_1^{\lf nt\rf}(1)\leq lu_1, 
\bigcap_{i=2}^{d}\Big\{\Zc^{\lf nt \rf}_i(1)\in (la_i, lb_i]\Big\},
\Lambda_{nt},\Xc_1(k)>0,\lfloor nt\rfloor<k\leq n\Big]\nonumber\\
 &=\frac{1}{\Po[\Lambda_n]}\int_0^{lu_1}\int_{la_2}^{lb_2}\cdots 
\int_{la_d}^{lb_d} \Po\Big[\Zc_1^{\lf nt\rf}(1)\in dx_1, \bigcap_{i=2}^{d}
\Big\{\Zc^{\lf nt \rf}_i(1)\in dx_i \Big\},\Lambda_{nt},\Xc_1(k)>0,
\lfloor nt\rfloor<k\leq n\Big]\nonumber\\
 &=\frac{\Po[\Lambda_{nt}]}{\Po[\Lambda_n]}\int_0^{lu_1}\int_{la_2}^{lb_2}
\cdots \int_{la_d}^{lb_d} \Po\Big[\Xc_1(k)>0,\lfloor nt\rfloor<k\leq n
\mid \Zc_1^{\lf nt\rf}(1)\in dx_1, 
\bigcap_{i=2}^{d}\Big\{\Zc^{\lf nt \rf}_i(1)\in dx_i \Big\}\Big]\nonumber\\
 &\phantom{*****************}\times\Po\Big[\Zc_1^{\lf nt\rf}(1)\in dx_1, 
\bigcap_{i=2}^{d}\Big\{\Zc^{\lf nt \rf}_i(1)\in dx_i \Big\}\mid\Lambda_{nt}\Big]
\nonumber\\
 &=\frac{\Po[\Lambda_{nt}]}{\Po[\Lambda_n]}\int_0^{lu_1}\int_{la_2}^{lb_2}
\cdots \int_{la_d}^{lb_d} 
\Po^{x\sqrt{\lfloor nt \rfloor}}\Big[\Zc^n_1(s)>0, 0\leq s \leq 1- n^{-1}\lf nt \rf \Big]\nonumber\\
 &\phantom{*****************}\times\Po\Big[\Zc_1^{\lf nt\rf}(1)\in dx_1, 
\bigcap_{i=2}^{d}\Big\{\Zc^{\lf nt \rf}_i(1)\in dx_i \Big\}\mid\Lambda_{nt}\Big].
\end{align}
By~(\ref{FI3}), (\ref{lamb1}), (\ref{POIT4}), and~(\ref{TIOFFF2})
we have $\IP$-a.s.
\begin{equation}
\label{Premlemrt}
 \lim_{n\to \infty}
\frac{\Po[\Lambda_{nt}]}{\Po[\Lambda_n]}=t^{-1/2}.
\end{equation}
Using Theorem \ref{Theouni} and Dini's theorem on 
uniform convergence of non-decreasing sequences of continuous functions, we obtain
\[
\lim_{n\to \infty}
\Po^{z\sqrt{\lfloor nt \rfloor}}\Big[\Zc_1^n(s)>0, 0\leq s \leq
1-n^{-1}\lfloor nt\rfloor \Big]={\tilde
N}\Big(z_1\Big(\frac{t}{1-t}\Big)^{1/2}\Big)
\]
uniformly in~$z$ on every compact set of the form $[0,K]\times [-K, K]^{d-1}$.
By Proposition~\ref{propcltord}, we have 
\[
\lim_{n \to \infty}
\Po\Big[\Zc_1^{\lf nt\rf}(1)\leq x_1, \bigcap_{i=2}^{d}\Big\{\Zc^{\lf nt \rf}_i(1)
\leq x_i \Big\}\mid\Lambda_{nt}\Big]=\exp\Big(-\frac{x_1^2}{2}\Big) 
\prod_{i=2}^{d}\int_{-\infty}^{x_i}\frac{e^{-\frac{v^2}{2}}}{\sqrt{2\pi}}dv.
\]
Now, applying Lemma~2.18 of~\cite{Igle} to~(\ref{findimdist4}), 
we obtain 
\begin{align*}
\lefteqn{\lim_{n \to \infty}\Po\Big[\Zc^n_1(n^{-1}\lfloor nt\rfloor)\leq u,
\bigcap_{i=2}^{d}\Big\{\Zc^n_i(n^{-1}\lfloor nt\rfloor)\in (a_i, b_i]\Big\}\mid 
\Lambda_n\Big]}\phantom{*********}\nonumber\\
&= \int_0^{u_1t^{-1/2}}\int_{a_2t^{-1/2}}^{b_2t^{-1/2}}\dots 
\int_{a_d t^{-1/2}}^{b_dt^{-1/2}}t^{-1/2}
{\tilde N}\Big(x_1\Big(\frac{t}{1-t}\Big)^{1/2}\Big)x_1e^{-\frac{x_1^2}{2}}
\prod_{i=2}^{d}\frac{e^{-\frac{x_i^2}{2}}}{\sqrt{2\pi}}dx_1\dots dx_d          .
\end{align*}
 Finally, make the change of variables $y=t^{1/2}x$ to obtain the
desired result.
 \qed

The final step in showing convergence of the f.d.d.'s is
\begin{prop}
 We have $\IP$-a.s., for all $k\geq 1$, $u_i>0$,  $-\infty<a_j^i< b_j^i<\infty$, 
$i\in[\![1,k]\!]$, $j\in[\![2,d]\!]$ and
$0<t_1<t_2<\dots<t_k\leq1$,
\begin{align}
\label{FDD2}
 \lefteqn{\lim_{n\to \infty} \Po\Big[\bigcap_{i=1}^{k}\Big\{\Zc^n_1(t_i)\leq
u_i,\Zc^n_2(t_i)\in (a^i_2, b^i_2],\dots,\Zc^n_d(t_i)\in (a^i_d, b^i_d]\Big \}\mid
\Lambda_n\Big]}\phantom{********}\nonumber\\
 &=\prod_{j=2}^{d}\int_{a^1_j}^{b^1_j}\dots \int_{a^k_j}^{b^k_j}
\frac{e^{-\frac{x_1^2}{2t_1}}}{\sqrt{2\pi t_1}} 
\frac{e^{-\frac{(x_2-x_1)^2}{2(t_2-t_1)}}}{\sqrt{2\pi (t_2-t_1)}}\dots 
\frac{e^{-\frac{(x_k-x_{k-1})^2}{2(t_k-t_{k-1})}}}
{\sqrt{2\pi (t_k-t_{k-1})}}dx_k\dots dx_1\nonumber\\
 &\phantom{**}\times \int_{0}^{u_1}
\dots\int_{0}^{u_k}q(0,0;t_1,x_1)q(t_1,x_1;t_2, y_2)\dots 
q(t_{k-1},x_{k-1};t_k,x_k)dx_k\dots dx_1.
\end{align}
\end{prop}

\noindent
\textit{Proof.}
The proof is by induction in~$k$. This result holds for $k=1$ by
virtue of~(\ref{FFD1}). Suppose~(\ref{FDD2}) is true for $k=m-1$, we
show that it can be extended to $k=m$. Let $t'_i=n^{-1}\lf t_i n\rf$ and let 
\[
\Dc_i=\{x\in \R^d: x_1\leq u_1, a^i_j<x_i\leq b^i_j, j\in [\![2,d]\!]\}
\]
for $i\in [\![1,m]\!]$. We mention here that in this proof, $y^i$ for 
$i\in[\![1,m]\!]$ are all elements of $\R^d$ while $y_i$ for $i\in[\![1,m]\!]$ 
belong to~$\R$. By 
the same argument as in the beginning of the proof of 
Proposition~\ref{prop_Z<x}, observe that 
\begin{equation}
\label{fraaa}
\lim_{n\to \infty}\Po\Big[\bigcap_{i=1}^{m}\{\Zc^n(t_i)\in \Dc_i\}\mid
\Lambda_n\Big]=\lim_{n\to \infty} \Po\Big[\bigcap_{i=1}^{m-2}\{\Zc^n(t_i)\in \Dc_i\},
\bigcap_{i=m-1}^{m}\{\Zc^n(t'_i)\in \Dc_i
 \}\mid
\Lambda_n\Big]
\end{equation}
provided that the limits exist.
Then, we write for sufficiently large~$n$
\begin{align}
\label{FDD3}
\lefteqn{
\Po\Big[\bigcap_{i=1}^{m-2}\{\Zc^n(t_i)\in \Dc_i\},\bigcap_{i=m-1}^{m}\{\Zc^n(t'_i)
\in \Dc_i  \}\mid
\Lambda_n\Big]
}\phantom{**}\nonumber\\
&=\frac{1}{\Po[\Lambda_n]} \int_{\Dc_{m-1}}\int_{\Dc_m}
\Po\Big[\Zc^n(t_1)\in  \Dc_1,\dots, \Zc^n(t_{m-2})\in  \Dc_{m-2},
\nonumber\\
 &\phantom{****************} \Zc^n(t'_{m-1}) \in
dy^{m-1}, \Zc^n(t'_m)\in
dy^{m}, X_1(1)>0,\dots,X_1(n)>0\Big]\nonumber\\
 &=\frac{\Po[\Lambda_{ nt_{m-1}}]}{\Po[\Lambda_n]}
\int_{\Dc_{m-1}}\int_{\Dc_m}
 \Po\Big[\Zc^n(t_1)\in \Dc_1,\dots, \Zc^n(t_{m-2})\in  \Dc_{m-2},
\Zc^n(t'_{m-1})\in dy^{m-1}\mid
\Lambda_{nt_{m-1}}\Big]\nonumber\\
 &\phantom{***************}\times \Po^{y^{m-1}\sqrt{n}}\Big[\Zc_1^n(s)>0, 0\leq s\leq t'_m-t'_{m-1},
\Zc^n(t'_m-t'_{m-1})\in
dy^m\Big] \nonumber\\
&\phantom{***************}\times
\Po^{y^m\sqrt{n}}\Big[\Zc_1^n(s)>0, 0\leq s\leq  1-t'_m\Big].
\end{align}
By the induction hypothesis we have
\begin{align}
\label{FDD4}
 \lefteqn{\lim_{n\to \infty} \Po\Big[\Zc^n(t_1)\in  \Dc_1,\dots, \Zc^n(t_{m-2})\in 
 \Dc_{m-2},
\Zc^n(t'_{m-1})\in \Dc_{m-1}\mid
\Lambda_{nt_{m-1}}\Big]}\phantom{************}\nonumber\\
 &=\prod_{j=2}^{d}\int_{a^1_j}^{b^1_j}\dots \int_{a^{m-1}_j}^{b^{m-1}_j}
\frac{e^{-\frac{y_1^2}{2t_1}}}{\sqrt{2\pi t_1}} 
\frac{e^{-\frac{(y_2-y_1)^2}{2(t_2-t_1)}}}{\sqrt{2\pi (t_2-t_1)}}\dots
 \frac{e^{-\frac{(y_{m-1}-y_{m-2})^2}{2(t_{m-1}-t_{m-2})}}}
{\sqrt{2\pi (t_{m-1}-t_{m-2})}}dy_{m-1}\dots dy_1\nonumber\\
&\phantom{**}\times \int_0^{u_1t_{m-1}^{-1/2}}\dots
\int_0^{u_{m-1}t_{m-1}^{-1/2}}q(0,0;t_1/t_{m-1},y_1)
q(t_1/t_{m-1},y_1;t_2/t_{m-1},y_2)\dots
\nonumber\\
 &\phantom{*****************}
q(t_{m-2}/t_{m-1},y_{m-2};1,y_{m-1})dy_{m-1}\dots dy_1.
\end{align}
On the other hand, by~(\ref{Premlemrt}) we have $\IP$-a.s.
\begin{equation}
\label{FDD5}
\lim_{n \to \infty} 
\frac{\Po[\Lambda_{nt_{m-1}}]}{\Po[\Lambda_n]}=t_{m-1}^{1/2}.
\end{equation}
Using Theorem~\ref{Theouni} and Dini's theorem on 
uniform convergence of non-decreasing sequences of continuous functions, we obtain
\begin{align}
\label{FDD6}
\lefteqn{
\lim_{n \to \infty}  \Po^{y^{m-1}\sqrt{n}}\Big[\Zc_1^n(s)>0, 0\leq s\leq t'_m-t'_{m-1},
\Zc^n(t'_m-t'_{m-1})\in \Dc_m\Big]}\phantom{****************}\nonumber\\
&= \prod_{j=2}^{d}\int_{a^m_j}^{b^m_j}\frac{e^{-\frac{(y_m-y^{m-1}_j)^2}
{2(t_m-t_{m-1})}}}{\sqrt{2\pi(t_m-t_{m-1})}}dy_m\times \int_0^{u_m} 
g(t_m-t_{m-1},y^{m-1}_1,v)dv
\end{align}
uniformly in $y^{m-1}$ on every compact set of the form $[0,K]\times [-K, K]^{d-1}$, and
\begin{equation}
\label{FDD7}
 \lim_{n \to \infty} \Po^{y^m\sqrt{n}}\Big[\Zc_1^n(s)>0, 0\leq s\leq  1-t'_m\Big]
={\tilde N}(y^m_1(1-t_m)^{-1/2})
\end{equation}
uniformly in~$y^m$ on every compact set of the form $[0,K]\times [-K, K]^{d-1}$.
Combining~(\ref{fraaa}), (\ref{FDD3}), (\ref{FDD4}), (\ref{FDD5}), 
(\ref{FDD6}), (\ref{FDD7}), and using Lemma~2.18 of~\cite{Igle}
twice, we obtain
\begin{align}
\label{FDD8}
 \lefteqn{
\lim_{n\to \infty}\Po\Big[\bigcap_{i=1}^{m}\{\Zc^n(t_i)\in \Dc_i\}\mid
\Lambda_n\Big]
}\phantom{***}\nonumber\\
 &= \prod_{j=2}^{d}\int_{a^1_j}^{b^1_j}\dots \int_{a^m_j}^{b^m_j}
\frac{e^{-\frac{x_1^2}{2t_1}}}{\sqrt{2\pi t_1}} 
\frac{e^{-\frac{(x_2-x_1)^2}{2(t_2-t_1)}}}{\sqrt{2\pi (t_2-t_1)}}\dots 
\frac{e^{-\frac{(x_m-x_{m-1})^2}{2(t_m-t_{m-1})}}}
{\sqrt{2\pi (t_m-t_{m-1})}}dx_m\dots dx_1\nonumber\\
&\phantom{**}\times  t_{m-1}^{-1}\int_0^{u_{m-1}}\int_0^{u_{m}}
\int_0^{u_1t_{m-1}^{-1/2}}\dots
\int_0^{u_{m-2}t_{m-1}^{-1/2}}q(0,0;t_1/t_{m-1},y_1)
q(t_1/t_{m-1},y_1;t_2/t_{m-1},y_2)\dots
\nonumber\\
 &\phantom{*************************}
q(t_{m-2}/t_{m-1},y_{m-2};1,y_{m-1}t^{-1/2}_{m-1}) dy_{m-1}\dots
dy_1\nonumber\\
 &\phantom{*************************} g(t_m-t_{m-1},y_{m-1}, y_m)
{\tilde N}(y_m(1-t_m)^{-1/2}) dy_{m}.
\end{align}
Now, make the change of variables  $t_{m-1}^{1/2}y_1=x_1,\dots,
t_{m-1}^{1/2}y_{m-2}=x_{m-2}$ in~(\ref{FDD8}) to obtain~(\ref{FDD2})
for $k=m$.
\qed


\subsection{Tightness}
\label{Tighteunesse}
 In this section, to finish the proof of~Theorem~\ref{Theocond},
 we prove that the sequence of measures $(\Po[\Zc^n\in \cdot \mid
\Lambda_n])_{n\geq 1}$ is tight $\IP$-a.s.
 First, we define the modulus of continuity for functions $f\in
C[0,1]$:
\[
w_f(\delta')=\sup_{|t-s|\leq \delta'}\{\|f(s)-f(t)\|_{\infty}\}
\]
where $s,t\in [0,1]$ and $\|\cdot\|_{\infty}$ is the $\infty$-norm 
on $\R^d$.
By Theorem~14.5 of~\cite{Kalen} it suffices
to show that $\IP$-a.s., for every $\hat{\eps}>0$
\begin{equation}
\label{gueuwint}
 \lim_{\delta' \downarrow 0}\limsup_{n\to
\infty}\Po[w_{\Zc^n}(\delta')\geq \hat{\eps} \mid \Lambda_n]=0
\end{equation}
since $\Zc^n(0)=0$.
Now observe that 
\begin{align}
\label{IOP}
\Po[w_{\Zc^n}(\delta')\geq \hat{\eps} \mid \Lambda_n]
 &=\Po\Big[\sup_{|t-s|\leq \delta'}\|\Zc^n(t)-\Zc^n(s)\|_{\infty}\geq \hat{\eps}\mid
\Lambda_n\Big]\nonumber\\
 &\leq\Po\Big[\sup_{|t-s|\leq 2\delta'}
\|\Xc( nt)-\Xc(ns)\|_{\infty}\geq \hat{\eps} \sqrt{n} \mid
\Lambda_n\Big]
\end{align}
for $n\geq 2/\delta'$.
 Let $m:=\lfloor 1/4\delta' \rfloor$ and divide the interval
$[0,1]$ into intervals 
$ I_k:=[\frac{k}{m}, \frac{k+1}{m}]$, for $0\leq
k\leq m-1$. Additionally, consider the intervals
$J_l:=[\frac{2l+1}{2m}, \frac{2l+3}{2m}]$, for $0\leq l\leq m-2$ and
$J_{m-1}:=\emptyset$. 
Observe that 
\begin{align}
\label{FLU1}
 \lefteqn{
\Po\Big[\sup_{|t-s|\leq 2\delta'}\|\Xc(nt)-\Xc(ns)\|_{\infty}\geq \hat{\eps} 
\sqrt{n} \mid
\Lambda_n\Big]}\phantom{**************}\nonumber\\
 &\leq \Po\Big[\Big\{\max_{k\leq m-1}\sup_{s,t \in { I}_k}
\|\Xc( nt)-\Xc(ns)\|_{\infty}\geq \hat{\eps}
\sqrt{n}\Big\}\nonumber\\
 &\phantom{*****}\cup  \Big\{\max_{l\leq m-1}\sup_{s,t \in J_l}
\|\Xc(nt)-\Xc(ns)\|_{\infty}\geq \hat{\eps} 
\sqrt{n}\Big\} \mid \Lambda_n\Big]\nonumber\\
 &\leq m \Big(\max_{ k\leq m-1} \Po\Big[\sup_{s,t \in { I}_k}
\|\Xc(nt)-\Xc(ns)\|_{\infty}\geq \hat{\eps} 
\sqrt{n}\mid \Lambda_n\Big]\nonumber\\
 &\phantom{*****}+ \max_{ l\leq m-1} \Po\Big[\sup_{s,t \in J_l}
\|\Xc(nt)-\Xc(ns)\|_{\infty}\geq \hat{\eps} 
\sqrt{n}\mid \Lambda_n\Big]\Big)
\end{align}
with the convention that $\sup_{s,t\in\emptyset}\{\cdot\}=0$.
Our next step is to bound from above the $\limsup_{n\to
\infty}$ of both terms in parentheses in the right-hand side 
of~(\ref{FLU1}). 
As an example, let us treat the terms indexed by $I_k$ for $k\in
[\![1,m-1]\!]$. The term indexed by~$I_0$ and those indexed by
$J_k$, $k\in [\![1,m-1]\!]$ can be treated in a similar way. To do
that, we will use the same approach as in
the proof of Proposition~\ref{propcltord}. Analogously to~(\ref{FI1}) we have 
for
$\eps\in (0,1)$ and $\delta\in (0,1)$,
\begin{align}
\label{cliff}
 \lefteqn{\Po\Big[\sup_{s,t \in { I}_k} \|\Xc(nt)-\Xc( ns)\|_{\infty} 
 \geq {\hat \eps}  \sqrt{n}\mid
\Lambda_n\Big]}\phantom{***********} \nonumber\\
& \leq (\Po[\Lambda_n])^{-1}\Po\Big[\sup_{s,t \in { I}_k} 
\|\Xc(nt)-\Xc( ns)\|_{\infty}  \geq {\hat \eps}  \sqrt{n},A_{0\to R},
\Lambda_n, \beta_{\Rc_{\eps, n}}\leq \delta nm^{-1}\Big]\nonumber\\
&\phantom{**}+\Po[\beta_{\Rc_{\eps, n}}> \delta nm^{-1}\mid \Lambda_n]
+\Po[\beta_{\Rc_{\eps, n}}> \eps^{1/2} n\mid \Lambda_n].
\end{align}
Analogously to~(\ref{FI2}), we obtain
\begin{align*}
 \lefteqn{
(\Po[\Lambda_n])^{-1}\Po\Big[\sup_{s,t \in { I}_k} \|\Xc(nt)-\Xc( ns)\|_{\infty} 
 \geq {\hat \eps}  \sqrt{n},A_{0\to R},\Lambda_n, \beta_{\Rc_{\eps, n}}
\leq \delta n m^{-1}\Big]
}\phantom{*************}\nonumber\\
 &\leq \frac{\Po[A_{0\to R}]}{\Po[\Lambda_n]}\max_{y\in R_{\eps, n}}
\max_{j\leq \lf \frac{\delta n}{m}\rf}\Po^y\Big[\sup_{s,t \in { I}_k} 
\|\Xc(nt-j)-\Xc( ns-j)\|_{\infty}  \geq {\hat \eps}  \sqrt{n}\Big].
\end{align*}
Now, observe that for all sufficiently large~$n$
\begin{align}
\max_{j\leq \lf \frac{\delta n}{m}\rf}\Po^y\Big[\sup_{s,t \in { I}_k} 
\|\Xc(nt-j)-\Xc( ns-j)\|_{\infty}  \geq {\hat \eps}  \sqrt{n}\Big]
 &\leq \Po^{y}\Big[\sup_{s,t \in I'_k} \|\Xc(nt)-\Xc(ns)\|_{\infty}
\geq {\hat \eps}  \sqrt{n}\Big]
\end{align}
with $I'_k=[\frac{k-2\delta}{m}, \frac{k+1}{m}]$.
 Now, let $I''_k=[\frac{k-3\delta}{m}, \frac{k+1}{m}]$. By
Theorem~\ref{Theouni} and the estimate on the tail of the
Gaussian law given in~\cite{PerMot}, Appendix~B, Lemma~12.9, we have
\begin{align}
\limsup_{n \to \infty}\max_{y\in R_{\eps,n}} \Po^{y}\Big[\sup_{s,t \in I'_k}
 \|\Xc(nt)-\Xc(ns)\|_{\infty}
\geq {\hat \eps}  \sqrt{n}\Big]
 &\leq d\cdot P\Big[\sup_{s,t \in I''_k} |W_1(t)-W_1(s)|\geq {\hat \eps}
\Big]\nonumber\\
 &\leq 8d\cdot P\Big[W_1\Big(\frac{1+3\delta}{m}\Big)\geq {\hat \eps}
\Big]\nonumber\\
 &\leq \frac{16d}{{\hat \eps} \sqrt{2\pi
m}}\exp\Big\{-\frac{{\hat \eps}^2m}{8}\Big\} 
\end{align}
since $\delta<1$.
We obtain
\begin{equation}
\label{RATA}
\limsup_{n \to \infty} \max_{y\in R_{\eps, n}}
\max_{j \leq \lf \frac{\delta n}{m}\rf}
\Po^y\Big[\sup_{s,t \in { I}_k} \|\Xc(nt-j)-\Xc( ns-j)\|_{\infty}  
\geq {\hat \eps}  \sqrt{n}\Big]
\leq \frac{16d}{\hat{\eps}\sqrt{2\pi m}}
\exp\Big\{-\frac{{\hat \eps}^2m}{8}\Big\}.
\end{equation}
Thus, we have by (\ref{FI3}), (\ref{URF}), (\ref{cliff}), and~(\ref{RATA})
\begin{align}
\label{FLU2}
 \lefteqn{
\limsup_{n \to \infty}\Po\Big[\sup_{s,t \in { I}_k}
\|\Xc( nt )-\Xc( ns )\|_{\infty}\geq {\hat \eps} 
 \sqrt{n}\mid \Lambda_n\Big]
}\phantom{***}\nonumber\\
 &\leq \Big(\frac{2\eps\sig_1}{\sqrt{2\pi}}+o(\eps)\Big)^{
-1} \Big(\frac{16d}{{\hat \eps} \sqrt{2\pi
m}}\exp\Big\{-\frac{{\hat \eps}^2m}{8}\Big\}\Big)+f(\eps)+g(\eps)
+ \limsup_{n \to \infty}\Po[\beta_{\Rc_{\eps, n}}> \delta nm^{-1}\mid \Lambda_n].
\end{align}
Combining~(\ref{FLU1}) and (\ref{FLU2}) we find
\begin{align}
\label{IOKL}
 \lefteqn{
\limsup_{n \to \infty}\Po\Big[\sup_{|t-s|\leq
2\delta'}\|\Xc(nt)-\Xc( ns )\|_{\infty}\geq {\hat
\eps}  \sqrt{n} \mid \Lambda_n\Big]
}\phantom{********************}\nonumber\\
 &\leq 2m\Big( \frac{16d}{{\hat \eps} \sqrt{2\pi
m}}\exp\Big\{-\frac{{\hat \eps}^2m}{8}\Big\} 
\Big(\frac{2\eps\sig_1}{\sqrt{2\pi}}+o(\eps)\Big)^{
-1}\nonumber\\
&\phantom{*****}+f(\eps)+g(\eps)
+ \limsup_{n \to \infty}\Po[\beta_{\Rc_{\eps, n}}> \delta nm^{-1}
\mid \Lambda_n]\Big).
\end{align}
Then, let $\eps=m^{-3}$ and $\delta=m^{-1/2}$ in~(\ref{IOKL}). 
We have by~(\ref{URF})
\begin{equation*}
\limsup_{n \to \infty}\Po[\beta_{\Rc_{\eps, n}}>
 \delta nm^{-1}\mid \Lambda_n]=\limsup_{n \to \infty}\Po[\beta_{\Rc_{\eps, n}}
> \eps^{1/2} n\mid \Lambda_n]\leq f(m^{-3})+g(m^{-3}).
\end{equation*}
Therefore, we obtain 
\begin{equation*}
 \lim_{m \to \infty}\limsup_{n \to \infty}
\Po\Big[\sup_{s,t \in {\hat I}_k}
\|\Xc( nt )-\Xc(ns)\|_\infty\geq {\hat \eps} 
\sqrt{n}\mid \Lambda_n\Big]=0.
\end{equation*}
 As ${\hat \eps}$ is arbitrary and $m=\lfloor 1/4\delta'\rfloor$, 
using~(\ref{IOP}), this last expression proves~(\ref{gueuwint}) and
consequently the tightness of the sequence $\big(\Po[\Zc^n\in \cdot
\mid \Lambda_n]\big)_{n\geq 1}$.
\qed


\section*{Acknowledgements}
This work was mainly done when N.G.\ was visiting Brazil in the end
of 2010; this visit was supported by FAPESP (2010/16085--7).
C.G.\ is grateful to FAPESP (grant 2009/51139--3) for financial
support. S.P.\ and M.V.\ were partially supported by
CNPq (grant 300886/2008--0 and 301455/2009--0). 
C.G., S.P., and M.V.\ also thank CNPq (472431/2009--9) and FAPESP
(2009/52379--8) for financial support. We also gratefully acknowledge the 
John von Neumann guest professor program of Technische Universit\"at M\"unchen which supported a visit of S.P.\ and M.V. to Munich.


\begin{thebibliography}{19}


\bibitem{Bar} \textsc{M.T.~Barlow} (1995)  
St Flour Lecture Notes: Diffusions on Fractals.  \textit{Lect. Notes Math.} \textbf{1690},
1--121.

\bibitem{BD} \textsc{M.T.~Barlow, J.-D.~Deuschel} (2010)  
Invariance principle for the random conductance model with unbounded
conductances. \textit{Ann.\ Probab.} \textbf{38} (1),
234--276.

\bibitem{BB} \textsc{N.~Berger, M.~Biskup} (2007) 
Quenched invariance principle 
for simple random walk on percolation clusters.  
\textit{Probab.\ Theory Related Fields\/} \textbf{137}, 83--120.

\bibitem{Billing} \textsc{P.~Billingsley} (1968)
\textit{Convergence of Probability Measures} (1st ed.).
Wiley, New York.

\bibitem{Bis} \textsc{M.~Biskup} (2011)
Recent progress on the random conductance model. 
\textit{Prob.\ Surveys} \textbf{8}, 294--373.

\bibitem{BP} \textsc{M.~Biskup, T.M.~Prescott} (2007)
Functional CLT for random walk among bounded random conductances. 
\textit{Elect.\ J.\ Probab.} \textbf{12}, paper No.~49, 1323--1348.

\bibitem{Del}
\textsc{T.~Delmotte} (1999)
Parabolic Harnack inequality and estimates of Markov chains on graphs.
\textit{Rev.\ Mat.\ Iberoamericana} \textbf{15} (1), 181--232.

\bibitem{CPSV2} \textsc{F.~Comets, S.~Popov, G.M.~Sch\"utz,
M.~Vachkovskaia} (2010)
 Quenched invariance principle for Knudsen stochastic billiard
in random tube.
\textit{Ann.\ Probab.} \textbf{38} (3), 1019--1061.

\bibitem{CPSV3} \textsc{F.~Comets, S.~Popov, G.M.~Sch\"utz,
M.~Vachkovskaia} (2010)
Knudsen gas in a finite random tube: transport diffusion and first
passage properties.
\textit{J.\ Statist.\ Phys.} \textbf{140}, 948--984.

\bibitem{GP} \textsc{C.~Gallesco, S.~Popov} (2012)
Random walks among random conductances I: Uniform quenched CLT.
\textit{Elect.\ J.\ Probab.} \textbf{17}, paper No.~85, 1--22.

\bibitem{GP2} \textsc{C.~Gallesco, S.~Popov}
Random walks among random conductances II: Conditional quenched CLT.
\textit{Available on Arxiv: 1210.0591.}


\bibitem{Igle} \textsc{D.~Iglehart} (1974)
Functional central limit theorems for random walks conditioned 
to stay positive.
 \textit{Ann.\ Probab.} \textbf{2} (4), 608--619.

\bibitem{Kalen} \textsc{O.~Kallenberg} (1997)
\textit{Foundations of modern probability.}
Springer, New York.

\bibitem{M} \textsc{P.~Mathieu} (2008)
Quenched invariance principles for random walks with random
conductances.
\textit{J.\ Statist.\ Phys.} \textbf{130} (5), 1025--1046.

\bibitem{PerMot} \textsc{P.~M\"orters, Y.~Peres} (2010)
\textit{Brownian Motion.} Cambridge University Press.

\bibitem{RY} \textsc{D.~Revuz, M.~Yor} (1999)
\textit{Continuous Martingales and Brownian Motion}.
Springer, Berlin.

\bibitem{Rh10} \textsc{R.~Rhodes} (2010)
Stochastic homogenization of reflected stochastic differential equations.
\textit{Electr.\ J.\ Probab.}  \textbf{15}, 989--1023.



\end{thebibliography}
\end{document}